
%
%
%
%
\documentstyle{amsppt}
\magnification=1200
\def\A{{\bold A}}
\def\a{{\frak a}}

\def\c{{\frak c}}
\def\m{{\frak m}}

\def\dash{^\prime}

\def\bs{\backslash}

\def\C{{\Bbb C}}

\def\dash{^\prime}
\def\det{\text{det}}

\def\End{\text{End}}

\def\hidehrule#1#2{\kern-#1%
  \hrule height#1 depth#2 \kern-#2 }
\def\hidevrule#1#2{\kern-#1{\dimen0=#1
  \advance\dimen0 by#2\vrule width\dimen0}\kern-#2 }

\def\il{\int\limits_}

\def\makeblankbox#1#2#3#4{\hbox{\vbox{\hidehrule{#1}{#2}%
  \kern-#1 %
  \hbox to#3{\hidevrule{#1}{#2}%
    \raise#4\vbox to #1{}%
    \vtop to #1{}%
    \hfil\hidevrule{#2}{#1}}%
  \kern-#1\hidehrule{#2}{#1}}}}

\def\ord{\text{ord}}

\def\O{{\Cal O}}

\def\Q{{\Bbb Q}}

\def\sdots{\mathinner{\mskip1mu\raise1pt\vbox{\kern7pt\hbox{.}}
    \mskip2mu\raise4pt\hbox{.}\mskip2mu\raise7pt\hbox{.}\mskip1mu}}
\def\sgn{\text{sgn}}

\def\smudge{\hbox{\quad\vrule height6pt width9.7pt}\medbreak}
\def\t#1{{^\top\!#1}}

\catcode`\@=11
\newcount\sectno \sectno=0
\def\docname=#1{\def\docn@me{#1}}
\newcount{\ye@r}
\docname={}
\newcount{\ye@r}
\def\shortd@te{\ye@r=\year\advance\ye@r-1900
  \number\month--\number\day--\number\ye@r}
\def\toplin@{{\tenpoint\smc\shortd@te\hfil\docn@me\
  \hfil\ifnum\sectno>0 \qquad\number\sectno--\fi\number\pageno}}
\headline=\toplin@
\output={\output@}
\def\output@{\shipout\vbox{\vglue3pc%
 \iffirstpage@ \global\firstpage@false
  \pagebody \makefootline%
 \else \ifrunheads@ \makeheadline \pagebody
       \else \pagebody \makefootline \fi
 \fi}%
 \advancepageno \ifnum\outputpenalty>-\@MM\else\dosupereject\fi}
\catcode`\@=\active
\magnification=1200
\pagewidth{32pc}\pageheight{41pc}
\loadeusm
\loadeurm
\def\A{{\Cal A}}
\def\B{{\Cal B}}
\def\H{{\eusm H}}
\def\Wh{{\eusm W}\!{\eurm h}}
\def\Wa{{\eusm W}\!{\eurm a}}
\def\dash{^\prime}
\def\il{\int\limits_}
\def\Ind{\text{\rm Ind}}
\def\Iwahori{\Cal B_2}
\def\teetwo{T^{\text s}_2}
\def\O{{\Cal O}}
\def\smudge{\hbox{\quad\vrule height6pt width9.7pt}\medbreak}
\docname={\smc Explicit Formulas for the Waldspurger and Bessel Models}
\topmatter
\title Explicit Formulas for the Waldspurger and Bessel Models\endtitle
\author Daniel Bump, Solomon Friedberg, and Masaaki Furusawa\endauthor
\affil Stanford University\\University of California, Santa Cruz\\
Mathematical Sciences Research Institute\endaffil
\address \kern-1pt Department of Mathematics,
Stanford University, Stanford, CA
94305-2125 \endaddress
\email bump\@gauss.stanford.edu \endemail
\address \kern-1pt Department of Mathematics,
University~of~California~Santa~Cruz,
Santa Cruz, CA  95064\endaddress
\email friedbe\@cats.ucsc.edu \endemail
\address {\kern-1pt
Mathematical Sciences Research Institute, 1000 Centennial
Drive, Berkeley, CA 94720-5070}
\endaddress
\email furusawa\@msri.org \endemail
\thanks *Research supported in part by National Science
Foundation grants DMS~9023441 (Bump) and DMS 9123845 (Friedberg),
by the AMS Centennial Research Fellowship (Bump), and by the NSF
Postdoctoral Research Fellowship DMS 9206242 (Furusawa).
Research at MSRI is supported in part by NSF grant
DMS 9022140.
\endthanks
\endtopmatter

\document

In this paper we will study certain models of irreducible admissible
representations of the split special orthogonal group $SO(2n+1)$ over
a nonarchimedean local field.  If $n=1$, these models were considered
by Waldspurger \cite{Wa1,Wa2}, and arose in his profound studies of
the Shimura correspondence.  If $n=2$, they were considered by
Novodvorsky and Piatetski-Shapiro \cite{NP}, who called them {\it
Bessel models}, and for general $n$ they were studied by Novodvorsky
\cite{No}.  In the works cited, these authors established the
uniqueness of these models; in this paper we establish functional
equations and explicit formulas for them.  In general, these models
arise from a variety of Rankin-Selberg integrals (for example, those of
Andrianov~\cite{An}, Furusawa~\cite{Fu1}, and
Sugano \cite{Su}), and the results of
this paper will naturally have applications to the study of
L-functions. Moreover, these models arise in the study of the theta
correspondence between $SO(2n+1)$ and the double cover of $Sp(2n)$,
and they will therefore be of importance in generalizing the work of
Waldspurger (see Furusawa~\cite{Fu2}).

In the final Section, we present a global application of the explicit
formulas:  we consider
the Eisenstein series (6.1) on $SO(2n+1)$ formed with a cuspidal
automorphic representation $\pi$ on $GL(n)$, and we
show that its Bessel period (6.2) is essentially
a product of L-series
$$L\big(n(s-1/2)+1/2,\pi\big)\,L\big(n(s-1/2)+1/2,\pi\otimes\eta\big),$$
where $\eta$ is a quadratic character.
This result generalizes work for $n=2$ and base field $\Q$
of Mizumoto~\cite{Mi} and B\"ocherer~\cite{B\"o}, and puts it in
a representation-theoretic context.

This global application is closely
related to the results of Bump, Friedberg, and Hoffstein~\cite{BFH2}.
That paper computes the spherical Whittaker functions on the
metaplectic double cover of $Sp(2n)$.
(Whittaker models on that group are also unique.)  The computation
in \cite{BFH2} has the following consequence:
if one forms the metaplectic Eisenstein series on the double
cover of $GSp(2n)$ with a cuspidal automorphic representation
$\pi$ of $GL(n)$ (which is possible because the cover splits over
$GL(n)\subset Sp(2n)$), the Whittaker coefficients of this Eisenstein
series are quadratic twists of the standard L-function of $\pi$. The
close relation between these two computations is a reflection of the
following result
of Furusawa~\cite{Fu2}, generalizing the case $n=2$ in
Piatetski-Shapiro and Soudry~\cite{PS}: the (special) Bessel
coefficient of a cusp form on $SO(2n+1)$ essentially agrees with the
Whittaker coefficient of the theta lift on the double cover of
$Sp(2n)$. If instead of a cusp form one considers the Eisenstein
series (6.1), the theta correspondent on the metaplectic group is the
metaplectic Eisenstein series, and our calculation implies that this
result of Furusawa for cusp forms is true for these Eisenstein
series also.  (Our calculation of the Bessel
period is in fact direct and independent of \cite{Fu2}.)

These results should have an application to the nonvanishing
of L-functions under quadratic twists.
Namely, there are Rankin-Selberg integrals on the
double cover of $GSp(6)$ (\cite{BG})
and on $SO(7)$ (\cite{Gi}) unfolding to Dirichlet series
involving the Whittaker (resp.\ Bessel) periods
described above, that is, to Dirichlet series whose individual
coefficients are the quadratic twists of a standard $GL(3)$
L-series. (The two constructions give Dirichlet series whose
individual coefficients are Euler products which agree at almost
all places.)
Arguing as in~\cite{BFH1},
one should be able to show that an infinite number of these
quadratic twists are nonzero.
In fact, these integrals are the next members of a series
beginning with Siegel's calculation of the Mellin transform of a
metaplectic $GL(2)$ Eisenstein series and including
integrals of Hecke type on the double cover of $GSp(4)$
(due in a nonmetaplectic context to Novodvorsky; see~\cite{BFH1}) and
on $SO(5)$ (due to Maass).
The elucidation of this scenario owes much to discussions with Duke,
Ginzburg, Goldfeld, and Hoffstein. In particular, the verification
that our
results could be applied to the evaluation of (6.2) was first worked
out in conversation with Ginzburg.

\subheading{1. Notations and Statement of Results}
Let $F$ be a nonarchimedean local field of characteristic different
from 2.  Let $\O$ denote the ring of integers of $F$, $\varpi$ denote
a local uniformizer, $q$ denote the cardinality of the residue field
$\O/\varpi\O$, and $|~|_F$ denote the absolute value on $F$,
normalized so that $|\varpi|_F=q^{-1}$.

We begin by describing our results on the Waldspurger model.
Let $G_2=GL(2,F)$, and
let $T_2$ be a maximal torus in $G_2$. Then $T_2$ is the connected
component of the identity in a group of orthogonal similitudes of
degree two corresponding to some quadratic form.
If $(\pi_2,V_{\pi_2})$ is an irreducible admissible
representation of $G_2$, and if $\sigma:T_2(F)\to\C$ is a character,
then there exists at most one linear functional $W:V_{\pi_2}\to\C$
(up to scalar multiplication) such that
$$W(\pi_2(t)\,v)=\sigma(t)\,W(v)\tag1.1$$
for all $t\in T_2$ and $v\in V_{\pi_2}$. This is proved when
$\sigma=1$ in Waldspurger~\cite{Wa1}, Proposition~9$'$, and
the proof in the general case is identical (as pointed
out in~\cite{Wa2}, Lemme~8).
In order for such a functional to exist, since $T_2$
contains the center $Z_2$ of $G_2$, it is necessary that the restriction
of $\sigma$ to $Z_2$ match the central character of $\pi_2$.
We will call a functional satisfying (1.1) a {\it Waldspurger functional}.
The {\it Waldspurger model} for $\pi_2$ will be the space of all
functions of the form $g\mapsto W(\pi_2(g)v)$ with $v\in V_{\pi_2}$.

First let us consider the case where $T_2=T_2^{\text a}$ is nonsplit.
We will limit ourselves to the case where $T_2^{\text a}$ has the
form
$$\left\{\pmatrix x&y\\\epsilon y&x\endpmatrix|x^2-y^2\epsilon\ne0\right\},
\tag1.2$$
where $\epsilon\in\Cal O^\times$ is a nonsquare. In this case, let
$T_2^{\text a}(\Cal O)=T_2^{\text a}\cap GL(2,\Cal O)$.
We will further assume that $\sigma$ is
trivial on $T_2^{\text a}(\Cal O)$. Since $T_2^{\text a}$ is generated
by $T_2^{\text a}(\Cal O)$ and by the center of $G_2$, on which $\sigma$
is to agree with the central character of $\pi_2$, it follows that $\sigma$
is uniquely determined by these conditions. (We note that if $T_2^{\text a}$
and $\sigma$ are obtained by localizing global data then these conditions will
be satisfied locally almost everywhere at places where the global quadratic
form defining $T_2$ is nonsplit.)

Suppose that $\pi_2$ is in the unramified principal series. Then the
contragredient representation of $\pi_2$ is also spherical, and the
$GL(2,\O)$-fixed vector in the contragredient
representation is clearly $T_2^{\text a}(\O)$-invariant. Since it also
has the correct transformation property with respect to the center of
$G_2$, it is thus a Waldspurger functional.
Thus if $\phi$ is the spherical vector in $V_{\pi_2}$, the function
$g\mapsto W(\pi_2(g)\phi)$ is the spherical function for $\pi_2$, which
is given by the Macdonald formula (see~\cite{Ca1}).

One of our results will be a formula analogous to the Macdonald formula
for the {\it split} Waldspurger functional. Thus let $T_2=T_2^{\text s}$
be a split torus of $G_2$. Specifically, we may take $\teetwo$ to be
the group of diagonal matrices in $G_2$; also, let $B_2$ be the Borel
subgroup consisting of upper triangular matrices in $G_2$.

Let us construct a Waldspurger functional for the
unramified principal series.
Let $\xi_1$, $\xi_2$ be unramified quasicharacters of $F^\times$, and
let $\xi$ be the character of $B_2$ given by
$$\xi\left(\pmatrix a&b\\&d\endpmatrix\right)=
\xi_1(a)\xi_2(d).\tag1.3$$
Suppose $\pi_2=\Ind(\xi)$ (we also write $\pi_2=\Ind(\xi_1,\xi_2)$)
is the representation obtained by normalized induction from the
character $\xi$.  Thus $V_{\pi_2}$ consists of the complex-valued
locally constant functions $f$ on $G_2$ such that
$$f(bg)=\xi(b)\,\delta_{B_2}^{1/2}(b)\,f(g)\tag1.4$$
for all $b\in B_2$, $g\in G_2$, where $\delta_{B_2}$
is the modular character of $B_2$, and $\pi_2$ is
the right regular representation.

Let $\sigma$ be an unramified quasicharacter of $F^\times$.
Extend $\sigma$ to $\teetwo$ (we use the same letter) by the formula
$\sigma\left(\left(
\smallmatrix ab&\\&b\endsmallmatrix\right)\right)=\sigma(a)\,
\xi_1\xi_2(b).$
Then a Waldspurger functional on $V_{\pi_2}$ is defined as follows.
Suppose $f\in V_{\pi_2}$.  Let
$$\aligned
W(f)&=\il{Z_2\bs \teetwo}f\left(\pmatrix 0&1\\1&0\endpmatrix
\pmatrix 1&1\\0&1\endpmatrix t\right)\,\sigma^{-1}(t)\,d^\times t\\
&=\il{F^\times}f\left(\pmatrix 0&1\\1&0\endpmatrix
\pmatrix 1&1\\0&1\endpmatrix
\pmatrix a&0\\0&1\endpmatrix\right)\,\sigma^{-1}(a)\,
d^\times a.
\endaligned\tag1.5$$
The Haar measure on $F^\times$ is normalized so that
the measure of $\O^\times$ is 1.
Suppose that $\xi_i(\varpi)=\gamma_i$ for $i=1,2$, and that
$\sigma(\varpi)=\tau$.
As we shall show in Section 2 below, the integral (1.5)
is absolutely convergent in the region
$$\left|\gamma_1\tau^{-1}\right|<q^{1/2},\qquad
\left|\gamma_2^{-1}\tau\right|<q^{1/2}.\tag1.6$$
For these representations, (1.1) holds.

Let $\phi_\xi$ be the $K_2=GL(2,\O)$-fixed vector in $V_{\pi_2}$
such that $\phi_\xi(I_2)=1$.  Define a function
$W\!a_\xi\colon G_2\to \C$ by the formula
$$W\!a_\xi(g)=W(\pi_2(g)\phi_\xi).$$
This function is analogous to the Whittaker function obtained
from the standard Whittaker functional.

Our main result on the Waldspurger functional gives
the analytic continuation of the
function $W\!a=W\!a_\xi$ to the full space of unramified quasicharacters
$\xi_1$, $\xi_2$, and $\sigma$,
and an explicit formula for its value.  To describe this, note that
$$W\!a_\xi(t_2g\kappa_2)=\sigma(t_2)\,W\!a(g)\tag1.7$$
for all $t_2\in \teetwo$, $\kappa_2\in K_2$.
Hence it suffices to determine $W\!a_\xi$ on a set of coset representatives
for the double cosets $\teetwo\bs G_2/K_2$.  Using the Iwasawa
decomposition, it follows that a set of coset representatives
is given by the matrices
$$\eta_k=\pmatrix 1&1\\0&1\endpmatrix \pmatrix \varpi^k&0\\0&1\endpmatrix,$$
with $k\geq0$.
In Section 2 we shall prove

\proclaim{Theorem~1.1}  Suppose that (1.6) holds.  Let
$$\Wa_\xi(g)=\frac{(1-\gamma_1\tau^{-1}q^{-1/2})(1-\gamma_2^{-1}\tau
q^{-1/2})}{1-\gamma_1\gamma_2^{-1}q^{-1}}W\!a_\xi(g).$$
Then $\Wa_\xi$ is given by the formula
$$\multline
\Wa_\xi(\eta_k)=(1-q^{-1})^{-1}q^{-k/2}\times\\
\left[\gamma_1^k\,\frac{(1-\gamma_2\tau^{-1}q^{-1/2})
(1-\gamma_1^{-1}\tau q^{-1/2})}{1-\gamma_1^{-1}\gamma_2}+
\gamma_2^k\,\frac{(1-\gamma_1\tau^{-1}q^{-1/2})
(1-\gamma_2^{-1}\tau q^{-1/2})}{1-\gamma_1\gamma_2^{-1}}\right].
\endmultline
\tag1.8$$
In particular, the function $\Wa_\xi$, originally defined as an integral
when the inequalities (1.6) hold,
has holomorphic continuation to all $\gamma_1,\gamma_2,\tau\in\C^\times$,
and is invariant under the interchange of $\gamma_1$ and $\gamma_2$.
\endproclaim

\noindent
Note that if $\gamma_1=\gamma_2$, this must be interpreted correctly: both
the numerator and the denominator in (1.8) vanish, but their ratio
is analytic, so the formula still makes sense.

We also address the analytic continuation of the Waldspurger
functionals.  Let $\Lambda$ be the domain of
$(\gamma_1,\gamma_2,\tau)\in(\C^\times)^3$ such that
$\tau\gamma_1^{-1}\ne q^{1/2}$ and $\tau\gamma_2^{-1}\ne q^{-1/2}$.
We note that for fixed $\xi_1$ and $\xi_2$, the space $\Ind(\xi)$ may
be identified with the space $C^\infty\big((B_2\cap K_2)\backslash
K_2\big)$ of locally constant functions on $K_2$ which are left
invariant by $B_2\cap K_2$; indeed, such a function may be uniquely
extended to a function on $G_2$ satisfying (1.4), and every element of
$\Ind(\xi)$ arises uniquely from an element of $C^\infty\big((B_2\cap
K_2)\backslash K_2\big)$ in this way.  We shall also prove in
Section~2:

\proclaim{Theorem~1.2} Fix an element
$f\in C^\infty\big((B_2\cap K_2)\backslash K_2\big)$.  For
$(\gamma_1,\gamma_2,\tau)$ satisfying (1.6), $f$ may be extended uniquely
to an element of $\Ind(\xi)$, and the function $W(f)$ defined by (1.5)
may thus be regarded as an analytic function of three variables
$(\gamma_1,\gamma_2,\tau)$. This function has analytic continuation to
all of $\Lambda$, so the domain of definition of $W:\Ind(\xi)\to\C$ may be
extended to all $\xi$ for which the parameters are in $\Lambda$.
If $(\gamma_1,\gamma_2,\tau)\in\Lambda$, then $W$ defines a (possibly
zero) Waldspurger functional on $\Ind(\xi)$.
\endproclaim

\noindent
The analytic continuation assertion here can actually be deduced from
Theorem~1.1, since if $\xi$ is in general position then the representation
$\pi_2$ is irreducible, and hence every element of $\Ind(\xi)$ is a linear
combination of right translates of the spherical vector. Thus there is
some overlap between these two theorems.  In fact, however, we shall give a
direct proof of Theorem~1.2, independent from Theorem~1.1.

Next let us describe our results concerning the Bessel models.  For
$r\geq2$ let $SO(r,F)$ denote the split group of determinant one
orthogonal matrices $$SO(r,F)=\{g\in SL(r,F)~|~(gx,gy)=(x,y)
\quad\text{for all}\quad x,y\in F^{r}\},$$
where $(\;,\;)$ is the quadratic form
$$(x,y)=\sum_{i=1}^r x_i\,y_{r+1-i}.$$
Let $G=SO(2n+1,F)$, and let
$U$ be the subgroup of $G$ consisting of upper triangular unipotent
matrices whose center $3\times3$ block is the identity. An element of $U$
is of the form $u=(u_{ij})$ with $u_{ij}=0$ for $i,j=n,n+1,n+2$ and
$i\neq j$.  Let $\psi$ be a character of $F$ of conductor $\O$.
Given $S=(a,b,c)\in \O$ such that
$b^2+2ac\neq0$, define a character $\theta_S$ of $U$
by the formula
$$\theta_S(u)=\psi(u_{12}+u_{23}+\cdots+u_{n-2,n-1}
+au_{n-1,n}+bu_{n-1,n+1}+cu_{n-1,n+2}).\tag1.9$$
Let $T$ be the subgroup of $G$ consisting of the matrices of the form
$$\pmatrix I_{n-1}&&\\&g&\\&& I_{n-1}\endpmatrix$$
which, acting by conjugation, stabilize $\theta_S$.  This constrains
$g$ to lie in a suitable torus in $SO(3)$.  Thus $T$ is a torus in
$G$, which may be either split or nonsplit over $F$, depending on $S$.
Note that $T$ normalizes $U$, and hence $R:=TU$ is again a subgroup of
$G$.  Let $\lambda$ be a character of $T$, and extend $\theta_S$ to a
character of $R$ by $\theta_S(tu)=\lambda(t)\,\theta_S(u)$.

Let $\pi\colon G\to\End(V_\pi)$ be an admissible representation of
$G$.  Then a {\it Bessel functional} on $\pi$ is a linear functional
$B\colon V_\pi\to\C$ such that
$$B(\pi(tu)v)=\theta_S(tu)\,B(v),\tag1.10$$
for all $t\in T$, $u\in U$, and $v\in V_\pi$.  As shown by Novodvorsky
\cite{No}, if $\pi$ is irreducible then the dimension of the space of
such functionals is at most 1.  Note also that, in view of the
isomorphism of $SO(3,F)$ and $PGL(2,F)$, a Bessel functional on
$SO(3,F)$ may be identified with a Waldspurger functional.

One may similarly define a Bessel functional on
the larger group of orthogonal similitudes; however,
since this group is the direct
product of $G$ with its center, there is
no gain in generality by doing so.

Once again we shall consider this notion for the
unramified principal series.
Let $\chi_1,\cdots,\chi_n$ be
unramified quasicharacters of $F^\times$.
We shall consider the principal series representation of $G$
$$\pi=\Ind(\chi_1,\cdots,\chi_{n}).$$
In our notation, this will be the representation on space of locally
constant functions $\Psi$ on $G$ which satisfy
$$\Psi(bg)=\delta_B^{1/2}(b)\,\left(
\prod_{i=1}^{n}\chi_i(y_i )\right)
\cdot\Psi(g)\tag1.11$$
for all
$$b=\pmatrix y_1&*&*&*&*&*&*\\&\ddots&*&*&*&*&*\\
&&y_{n}&*&*&*&*\\&&&1&*&*&*\\&&&&y_n^{-1}&*&*\\
&&&&&\ddots&*\\&&&&&&y_1^{-1}\\\endpmatrix$$
in the standard Borel subgroup $B$ of $G$.  Here
$$\delta_B(b)=\left|\prod_{i=1}^ny_i^{2n-2i+1}\right|_F$$
is the modular character of $B_F$. The group action is by right
translation.  Let us write $\alpha_i=\chi_i(\varpi)$.
If the $\alpha_i$ are in general position
($\alpha_i\neq \alpha_j^{\pm1}$ for $i\neq j$ and $\alpha_i^2\neq1$),
then the isomorphism class of this representation is
invariant under permutations of the $\alpha_i$, as well as transformations of
the form $\alpha_i\mapsto\alpha_i^{\pm1}$.
Also, let us write $K$ for the standard maximal
compact subgroup $SO(2n+1,\O)$ of $G$.

We may construct a Bessel functional explicitly as follows.  There are
two cases: $T$ nonsplit, and $T$ split.  In the split case, we shall
take $\lambda$ unramified, i.e.\ identically one on $T\cap K$.

Suppose first that $T$ is a nonsplit torus. We assume that $a\in
\O^{\times}$.  Write $T(\O)$ for the subgroup $T\cap K$.  For a
permutation $s$ in the symmetric group $S_{2n+1}$, we shall also use
$s$ to denote the corresponding signed permutation matrix in
$SL(2n+1,F)$ (the matrix with $\sgn(s)$ in the $(s(i),i)$ position and
$0$ elsewhere).  Let $w_1=(1,2n+1)(2,2n)\cdots(n-1,n+3)$.  Suppose
$\Psi\in V_\pi$.  Then we define
$$B(\Psi)=\il{T(\O)}\il{U}\Psi(w_1ut)\,\theta_S(u)^{-1}\,du\,dt,\tag1.12$$
and this is a Bessel functional on $V_\pi$. The integral is absolutely
convergent if the quasicharacters $\chi_i$ are in a suitable region.
Indeed, since $T(\O)$ is compact, comparing with the standard
intertwining operator $T_{w_1}$ defined in Section 3 below (see Lemma
3.1), one finds that if
$$|\alpha_1|<\cdots<|\alpha_{n-1}|<\min(|\alpha_n|,|\alpha_n^{-1}|),
\tag1.13$$
then the integral (1.12) is absolutely convergent.

For the second case, suppose that $T$ is instead split.
We may suppose that $S=(0,1,0)$.
Let us introduce the following notation.
For $x\in F$, let us write $n(x)$ for the unipotent matrix
$$n(x)=\pmatrix I_{n-1}&&&&\\&1&x&-x^2/2&\\&&1&-x&\\
&&&1&\\&&&&I_{n-1}\endpmatrix,$$
and for $a\in F^\times$, let us write $t(a)$ for the diagonal matrix
in $T$ given by
$$t(a)=\pmatrix I_{n-1}&&&&\\&a&&&\\&&1&&\\&&&a^{-1}&\\&&&&I_{n-1}
\endpmatrix.$$
The matrices $t(a)$ give all of $T$.  Also we
write $\lambda(a)=\lambda(t(a))$, $\beta=\lambda(\varpi)$.
Let $w_0=w_1(n,n+2)$ be (a representative for)
the long element of the Weyl group $\Omega$.
Then for $\pi$ as above, and $\Psi\in V_\pi$, we define
$$B(\Psi)=\il{F^\times}\il{U}\Psi(w_0\,n(1)\,u\,t(a))\,\theta_S(u)^{-1}
\,\lambda^{-1}(a)\,du\,d^\times a.\tag1.14$$
Then this is a Bessel functional on $V_\pi$.
Once again,
the Haar measure on $F^\times$ is normalized so that
the measure of $\O^\times$ is 1.
As we shall show in Section 3 below, the integral (1.14) is
absolutely convergent if
$$|\alpha_1|<\cdots<|\alpha_{n-1}|<\min(|\alpha_n|,|\alpha_n^{-1}|),
\qquad |\alpha_n|<q^{1/2}\min(|\beta|,|\beta^{-1}|).
\tag1.15$$

Our first pair of results concerns the functional equations satisfied
by $B$.  We shall show that,
in both cases, the Bessel functional may be extended to all characters
$\chi=(\chi_1,\cdots,\chi_n)$ by a
variation of the familiar process whereby the standard intertwining
operators are analytically continued, and has a functional equation under
certain transformations of these characters $\chi$.
More precisely, the Weyl group $\Omega$ of $G$ acts on the characters
$\chi=(\chi_1,\cdots,\chi_{n})$, or what is the same thing in the
unramified case, on the parameters $\alpha_1,\cdots,\alpha_{n}$.
In terms of these parameters, $\Omega$ is the group of
transformations of $(\alpha_1,\cdots,\alpha_n)\in(\C^\times)^{n}$ generated
by
$$(\alpha_1,\ldots,\alpha_{n-1},\alpha_{n}) \mapsto
(\alpha_1,\ldots,\alpha_{n-1},\alpha_{n}^{-1}),$$
and by the action of the symmetric group $S_{n}$ on
$(\alpha_1,\cdots,\alpha_{n})$.
The cardinality of $\Omega$ is $2^n n!$.

Let $\Phi_\chi\in V_\pi$ be the standard nonramified vector. This is
the unique function in $V_\pi$ taking value one on $K$.
We define a function
$H=H_{\chi}\colon G\to\C$ by
$$H(g)=B(\pi(g)\,\Phi_\chi).$$
This function is once again analogous to the Whittaker function obtained
from the standard Whittaker functional.

\proclaim{Theorem~1.3}  Suppose that $T$ is nonsplit.
Then the function $H_{\chi}$, originally defined as an integral when
$$|\alpha_1|<\cdots<|\alpha_{n-1}|<
\min(|\alpha_{n}|,|\alpha_{n}^{-1}|),$$
has a meromorphic continuation to all
nonzero complex $\alpha_1$,\dots, $\alpha_{n}$.  Moreover the function
$$\H_{\chi}(g)=
\prod_{1\leq i<j\leq n}(1-\alpha_i\alpha_jq^{-1})^{-1}
(1-\alpha_i\alpha_j^{-1}q^{-1})^{-1}\,H_{\chi}(g)$$
is invariant under the action of $\Omega$ on the $\alpha_i$, and
is holomorphic for all $(\alpha_1,\dots,\alpha_n)\in(\C^\times)^{n}$.
\endproclaim

\proclaim{Theorem~1.4}  Suppose that $T$ is split.
Then the function $H_{\chi}$, originally
defined as an integral when
$$|\alpha_1|<\cdots<|\alpha_{n-1}|<\min(|\alpha_n|,|\alpha_n^{-1}|),
\qquad |\alpha_n|<q^{1/2}\min(|\beta|,|\beta^{-1}|),$$
has a meromorphic continuation to all
nonzero complex $\alpha_1$,\dots, $\alpha_{n}$, $\beta$.
Moreover the function
$$\H_{\chi}(g)=
{{\prod_{i=1}^n(1-\alpha_i\beta q^{-1/2})(1-\alpha_i\beta^{-1} q^{-1/2})}
\over{
\prod_{1\leq i<j\leq n}(1-\alpha_i\alpha_jq^{-1})
(1-\alpha_i\alpha_j^{-1}q^{-1})\,\prod_{i=1}^n(1-\alpha_i^2 q^{-1})}}
\,H_{\chi}(g)$$
is invariant under the action of $\Omega$ on the $\alpha_i$, and
is holomorphic for all $(\alpha_1,\dots,\alpha_n)\in(\C^\times)^{n}$
and $\beta$ satisfying $q^{-1/2}<\min(|\beta|,|\beta|^{-1})$.
\endproclaim

\noindent
The proofs of Theorems 1.3 and 1.4 are given in Section 3 below.

Our next pair of results gives an explicit formula for $\H_{\chi}(g)$.
As in the case of the Waldspurger model (see (1.7)), since
$$\H_\chi(rg\kappa)=\theta_S(r)\H_\chi(g)\tag1.16$$
for all $r\in R$, $\kappa\in K$, it suffices to determine $\H$ on a
set of coset representatives for $R\bs G/K$.  Let $k=(k_1,\cdots,
k_n)$ be a vector of integers with $k_1\geq0$.  If $T$ is nonsplit
then one finds, using the Iwasawa and Cartan decompositions (compare
Sugano \cite{Su}, Lemma 2-4), that a set of coset representatives for
$R\bs G/K$ is given by the diagonal matrices of the form
$$ d_k=\text{diag}(\varpi^{k_n'},
\varpi^{k_{n-1}'}, \cdots,\varpi^{k_1'},1,
\varpi^{-k_1'},\cdots,\varpi^{-k_{n-1}'},
\varpi^{-k_n'}),
$$
with
$$k_i'=k_{1}+\cdots+k_{i}.$$
If $T$ is split, then one finds instead that a set of coset
representatives is given by the matrices $g_k$ of the form
$$g_k=\cases d_k&\text{if $k_1=0$}\\
n(1)\,d_k&\text{if $k_1>0$.}\endcases
$$
For convenience, we write $h(k_1,\cdots,k_{n})$ for the
quantity $\H(d_k)$ (resp.\ $\H(g_k)$).

It follows from equation (1.16)
$H(g)=0$ unless $\theta_S$ is identically one on $R\cap gKg^{-1}$.
A short calculation shows that this condition implies that
$h(k_1,\cdots,k_{n})=0$ unless each $k_i\geq0$.

Let $\A$ be the alternator $\sum_{w\in\Omega} (-1)^{\text{length}(w)}\,w$
in the group algebra $\C[\Omega]$. Let
$\Delta=(-1)^n\A(\alpha_1^{n}\alpha_2^{n-1}\cdots\alpha_{n})$.
According to Weyl's identity for $Sp(2n,\C)$
$$\Delta=\prod_{i=1}^{n}\alpha_i^{-1+i-n}(1-\alpha_i^2)
\prod_{1\leq i<j\leq n}(1-\alpha_i\alpha_j)(1-\alpha_i\alpha_j^{-1})
.\tag1.17$$
Also let
$$e_k=-\frac{1}{2}\sum_{i=1}^n\left(n^2-(i-1)^2\right)\,k_{i}.$$
Then the evaluation of $h(k_1,\cdots,k_{n})$ is given by

\proclaim{Theorem~1.5} Suppose $T$ is nonsplit and $k_i\geq0$ for
$i=1$ to $n$.  Then
$$h(k_1,\cdots,k_{n})=
q^{e_k}(1+q^{-1})^{-1}\,\Delta^{-1}\,
\A\!\left(\,\,\prod_{i=1}^{n}
\alpha_{n+1-i}^{-k'_i-i}(1-\alpha_i^{2}q^{-1})
\right).$$
In particular, $\H_\chi\left(I_{2n+1}\right)=1$.
\endproclaim

\proclaim{Theorem~1.6} Suppose $T$ is split and $k_i\geq0$ for
$i=1$ to $n$.
\roster
\item
If $k_1=0$ then
$$\multline h(0,k_2,\cdots,k_{n})=
q^{e_k}\,\Delta^{-1}\\\times
\A\!\left(\alpha_n^{-1}
\prod_{i=1}^{n-1}
\alpha_{i}^{-k'_{n+1-i}-n-1+i}(1-\alpha_i\beta q^{-1/2})
(1-\alpha_i\beta^{-1}q^{-1/2})
\right).\endmultline$$
\item
If $k_1>0$ then
$$\multline h(k_1,k_2,\cdots,k_{n})=
q^{e_k}(1-q^{-1})^{-1}\,\Delta^{-1}\\\times
\A\!\left(\prod_{i=1}^{n}
\alpha_{n+1-i}^{-k'_i-i}(1-\alpha_i\beta q^{-1/2})
(1-\alpha_i\beta^{-1}q^{-1/2})
\right).\endmultline$$
\endroster
In particular, the function $\H_\chi$ is holomorphic for all
$(\alpha_1,\dots,\alpha_n,\beta)\in(\C^\times)^{n+1}$, and
$\H_\chi\left(I_{2n+1}\right)=1$.
\endproclaim
These theorems are proved in Section 4 below.

Let us finally formulate the meromorphic continuation of the Bessel
functional to all values of $\chi$. We will
formulate this result only in the split case; the nonsplit case is
nearly identical.  Suppose that
$\Psi\in C^\infty\big((B\cap K)\bs K\big)$. Given $\chi$, there is a unique
extension $\Psi_\chi$ of $\Psi$ to an element of $\Ind(\chi)$
satisfying (1.11).

\proclaim{Theorem~1.7} There exists a dense open subset $\Gamma$ of
$(\C^\times)^{n+1}$ such that there exists a Bessel functional $B$ on
$V_\pi$ for all
$(\alpha_1,\cdots,\alpha_n,\beta)\in\Gamma$; if (1.15) is satisfied,
this functional agrees with that defined by (1.14); and if $\Psi\in
C^\infty\big((B\cap K)\bs K\big)$, then $B(\Psi_\chi)$ is a
meromorphic function of $\alpha_1,\cdots,\alpha_n,\beta$, whose polar
set is contained in the complement of $\Gamma$.
\endproclaim

\noindent
As with Theorem~1.2, there is some overlap between this result and our
previous Theorems. We will prove Theorem~1.7 by means of a theorem of
Bernstein \cite{Be} in Section~5. Actually Bernstein's theorem implies
that the complement of $\Gamma$ may be taken to be a countable union of
hyperplanes.

\bigbreak
In view of the isomorphism of $SO(5,F)$ to the projective group of
symplectic similitudes $PGSp(4,F)$, the theorems above imply the
functional equation and the explicit formula for the Bessel model
on $G_4:=GSp(4)$. The Bessel model on $GSp(4,F)$ was first studied
by Novodvorsky and Piatetski-Shapiro \cite{NP}.
A formula for the generating series of this
model was given by Sugano in \cite{Su}, Proposition 2-5.
To state the explicit formula, consider the group
$$G_4=\left\{g\in GL(4,F)~\big|~ ^tg\pmatrix 0&I_2\\-I_2&0\endpmatrix
g=\nu(g)\pmatrix 0&I_2\\-I_2&0\endpmatrix, \nu(g)\in F^\times\right\}.$$
Let $S\in M_2(\O)$ be a nonsingular
symmetric matrix $S=\left(\smallmatrix c&b/2\\
b/2&a\endsmallmatrix\right)$.
Let $U_4$ be the unipotent radical of the Siegel parabolic
$$U_4=\left\{u=\pmatrix I_2&X\\0&I_2\endpmatrix~\big|~ X=^t\!X\right\},$$
and let $\theta_S$ be the character of $U_4$ given by
$\theta_S(u)=\psi(\text{tr}(SX))$.
Let $T$ be the torus in $G_4$ consisting of those matrices of the
form
$$\left(\matrix h&0\\0&\det(h)\,^th^{-1}\endmatrix\right)$$
which, acting by conjugation, stabilize $\theta_S$; thus $h\in G_2$
is required to satisfy $^th S h=\det(h) S$.
Then $T$ normalizes $U$, and hence $R:=TU$
is a subgroup of $G_4$.
Let $\lambda$ be a character of $T$,
and once again extend $\theta_S$ to a character of
$R$ by $\theta_S(tu)=\lambda(t)\,\theta_S(u)$.

If $\pi$ is an admissible representation of $G_4$,
a Bessel functional on $\pi$ is once again a linear
functional $B\colon V_\pi\to\C$ such that
$$B(\pi(tu)v)=\theta_S(tu)\,B(v),$$
for all $t\in T$, $u\in U$, and $v\in V_\pi$.
If $\pi$ transforms by a central character, then since
$T$ contains the center $Z_4$ of $G_4$, this notion
requires that $\lambda\mid_{Z_4}$ match this character.

Let $\pi=\Ind(\chi)$ be the principal series representation
obtained by normalized induction from a character $\chi$
of the standard Borel subgroup of $G_4$.  For $\chi$ in
a suitable domain a Bessel functional may be obtained
by integration as above.  If $\Phi_\chi\in V_\pi$ is the
standard nonramified vector, then define $H=H_\chi\colon
G_4\to\C$ by the formula $H(g)=B(\pi(g)\Phi_\chi)$.

To describe the explicit formula for this function,
define parameters $\alpha_1$ to $\alpha_4$ by
$$\align
\chi\pmatrix \varpi&&&\\&\varpi&&\\&&1&\\&&&1\endpmatrix=\alpha_1 \qquad
&\chi\pmatrix \varpi&&&\\&1&&\\&&1&\\&&&\varpi\endpmatrix=\alpha_2 \\
\chi\pmatrix 1&&&\\&1&&\\&&\varpi&\\&&&\varpi\endpmatrix=\alpha_3 \qquad
&\chi\pmatrix 1&&&\\&\varpi&&\\&&\varpi&\\&&&1\endpmatrix=\alpha_4.
\endalign$$
The Weyl group $\Omega_4$ acts on these parameters through all
permutations of the $\alpha_i$ which preserve the relation
$\alpha_1\alpha_3=\alpha_2\alpha_4$.
Let $\A_4$ be the alternator $\sum_{w\in\Omega_4} (-1)^{\text{length}(w)}\,w$
in the group algebra $\C[\Omega_4]$.
Observe moreover that the function  $H$ is completely determined
by its values on the elements
$$a_{k,l}=\pmatrix\varpi^{k+2l}&&&\\&\varpi^{k+l}&&\\&&1&\\&&&\varpi^{l}
\endpmatrix$$
with $k,l\geq0$ if $T$ is nonsplit, and by its values on
the elements
$$b_{k,l}=\cases a_{k,0}&\text{if $l=0$}\\m(1)a_{k,l}&\text{otherwise}
\endcases$$
with $k,l\geq0$ if $T$ is split, where we set
$$m(1)=\pmatrix 1&1&&\\&1&&\\&&1&\\&&-1&1\endpmatrix.$$
Then we have

\proclaim{Corollary 1.8} Suppose that $T$ is nonsplit and $a\in\O^\times$.
\roster
\item
The function $H_{\chi}$
may be defined by continuation for all unramified characters
$\chi$.  Moreover the function
$$\H_{\chi}(g)=
(1-\alpha_1\alpha_2^{-1}q^{-1})^{-1}(1-\alpha_2\alpha_3^{-1}q^{-1})^{-1}
H_\chi(g)$$
is holomorphic and invariant under the action of $\Omega$.
\item  The function $\H_\chi$ is given by the explicit formula
$$\H_\chi(a_{k,l})=
(1+q^{-1})^{-1}q^{-3k/2-2l}
\frac{\A_4\left(\alpha_3^{k+l+2}\alpha_2^l\alpha_4^{-1}
(1-\alpha_1\alpha_3^{-1}q^{-1})(1-\alpha_4\alpha_2^{-1}q^{-1})
\right)}{\A_4\left(\alpha_3^2\alpha_4^{-1}\right)}$$
valid for $k,l\geq0$.
\endroster
\endproclaim

In the split case, let us suppose without loss that $a=c=0$, $b=1$.
Let $\beta_1$, $\beta_2$ be the parameters
$$\beta_1=\lambda\left(
\pmatrix \varpi&&&\\&1&&\\&&1&\\&&&\varpi\endpmatrix\right),\qquad
\beta_2=\lambda\left(
\pmatrix 1&&&\\&\varpi&&\\&&\varpi&\\&&&1\endpmatrix\right).
$$
We are concerned with characters $\chi$ such that
$\alpha_1\alpha_3=\beta_1\beta_2$.
Then we have
\proclaim{Corollary 1.9} Suppose that $T$ is split.
\roster
\item
The function $H_{\chi}$
may be defined by continuation for all unramified characters
$\chi$, $\lambda$ such that $\chi\mid_{Z_4}=\lambda\mid_{Z_4}$.
Moreover the function
$$\H_{\chi}(g)=
{{\prod_{i,j=1}^2(1-\alpha_i\beta_j^{-1} q^{-1/2})}
\over{
(1-\alpha_1\alpha_2^{-1}q^{-1})
(1-\alpha_1\alpha_3^{-1}q^{-1})
(1-\alpha_2\alpha_3^{-1}q^{-1})
(1-\alpha_2\alpha_4^{-1}q^{-1})}}
\,H_{\chi}(g)$$
is holomorphic and invariant under the action of $\Omega$.
\item  The function $\H_\chi$ is given by the explicit formula
$$\align
\H_\chi(b_{k,0})=&
q^{-3k/2}
\frac{\A_4\left(\alpha_3^{k+2}\alpha_4^{-1}
\prod_{j=1}^2(1-\alpha_1\beta_j^{-1}q^{-1/2})\right)}
{\A_4\left(\alpha_3^2\alpha_4^{-1}\right)}\\
\intertext{valid for $k\geq0$ and by}
\H_\chi(b_{k,l})=&
(1-q^{-1})^{-1}\,q^{-3k/2-2l}
\frac{\A_4\left(\alpha_3^{k+l+2}\alpha_2^l\alpha_4^{-1}
\prod_{i=1,4; j=1,2}(1-\alpha_i\beta_j^{-1}q^{-1/2})\right)}
{\A_4\left(\alpha_3^2\alpha_4^{-1}\right)}\\
\endalign$$
valid for $k\geq0$, $l>0$.
\endroster
\endproclaim

\subheading{2. The Waldspurger Functional}
We begin the study of the split Waldspurger functional with
the following Lemma.

\proclaim{Lemma~2.1}  Suppose that the inequalities
(1.6) hold.  Then the integral (1.5) is absolutely convergent.
\endproclaim

\demo{Proof} Let $f\in\Ind(\xi)$ be given.
Since $K_2$ is compact, there is a number $C$ such that
$|f(\kappa)|\leq C$ for all $\kappa\in K_2$. One has
the matrix identities
$$\pmatrix 0&1\\1&0\endpmatrix
\pmatrix 1&1\\0&1\endpmatrix
\pmatrix a&0\\0&1\endpmatrix
=\cases \pmatrix a&1\\0&1\endpmatrix\pmatrix -1&0\\a&1\endpmatrix&
\text{if $|a|_F\leq 1$}\\
\pmatrix1&0\\0&a\endpmatrix\pmatrix0&1\\1&a^{-1}\endpmatrix&
\text{if $|a|_F>1$,}\endcases\tag2.1$$
where the last matrix in each case is in $K_2$.
Since $f\in\Ind(\xi)$, one obtains the inequalities
$$\left|f\left(\pmatrix 0&1\\1&0\endpmatrix
\pmatrix 1&1\\0&1\endpmatrix
\pmatrix a&0\\0&1\endpmatrix\right)\right|\leq
\cases C|a|_F^{1/2}|\xi_1(a)|&\text{if $|a|_F\leq1$}\\
C|a|_F^{-1/2}|\xi_2(a)|&\text{if $|a|_F>1$.}\endcases$$
Thus one sees that the integral (1.5) is majorized by the sum of the
integrals
$$\il{|a|_F\leq1}|\xi_1(a)|\,|\sigma^{-1}(a)|\,|a|_F^{1/2}\,d^\times a$$
and
$$\il{|a|_F>1}|\xi_2(a)|\,|\sigma^{-1}(a)|\,|a|_F^{-1/2}\,d^\times a.$$
But when the inequalities (1.6) hold, each of these integrals is an
absolutely convergent geometric series.  The Lemma follows.
\smudge
\enddemo

We turn to the proof of Theorem 1.1.
The proof makes essential use of ideas of Casselman and
Shalika~\cite{CS}, and of Banks~\cite{Ba}.

Let us define certain elements of the
representation $\pi=\pi_2=\text{Ind}(\xi)$ as follows. Let $\phi_\xi$ be
the normalized $K_2$-fixed vector, and for $k\ge0$, let
$$F_k(g)=\int_{\Cal O}\phi_\xi\!\left(g\pmatrix1&z\\&1\endpmatrix
\pmatrix\varpi^k\\&1\endpmatrix\right)dz,\tag2.2$$
$$\zeta_k(g)=\phi_\xi\!\left(g\pmatrix 1&1\\&1\endpmatrix
\pmatrix\varpi^k\\&1\endpmatrix\right).$$
It is easy to see that in the integral (2.2), the integrand has constant
value when $z$ has constant valuation; and if $z$ is a unit, this constant
value is $\zeta_k(g)$. Thus
$$(1-q^{-1})\zeta_k(g)=F_k(g)-
\int_{\varpi\Cal O}\phi_\xi\!\left(g\pmatrix1&z\\&1\endpmatrix
\pmatrix\varpi^k\\&1\endpmatrix\right)dz.$$
A simple change of variables shows that the second integral equals
$$q^{-1}\,\pi\!\pmatrix\varpi\\&1\endpmatrix F_{k-1}(g),$$
and so applying the Waldspurger functional $W$, we obtain
$$(1-q^{-1})W(\zeta_k)=W(F_k)-q^{-1}\tau W(F_{k-1}).\tag2.3$$

As in Casselman~\cite{Ca1} and Casselman and Shalika~\cite{CS}, the
vectors $F_k$ are fixed by the Iwahori subgroup $\Iwahori$ of $G_2$. We
remind the reader of the Casselman basis of the Iwahori fixed vectors of
$G_2$. Assuming that $\xi$ is regular, we define linear functionals
$T_w$ for $w$ a Weyl group representative of $G_2$ by
$$T_wf=\int_{N\cap wNw^{-1}\bs N}f(w^{-1}n)\,dn,$$
where $N$ is the group of upper triangular unipotent matrices in $G_2$.
These functionals are linearly independent on the Iwahori fixed vectors,
and the Casselman basis $f_w$ is defined by $T_wf_{w'}=\delta(w,w')$
(Kronecker delta) for $w$ and $w'$ in the Weyl group. If $w_0$ is (a
representative of) the long element of the Weyl group, then
$$f_{w_0}(g)=\cases\phi_\xi(g)&\text{if $g\in B_2w_0\Iwahori$,}\\
0&\text{otherwise.}\endcases\tag2.4$$
The element $f_1$ is given by a more complicated formula, and we do not
need to know it. Since $F_k$ is an Iwahori fixed vector, we can write
$$F_k=c(1,\xi)\,f_1+c(w_0,\xi)\,f_{w_0},$$
and by definition of the $f_w$, $c(w,\xi)=T_wF_k$. It is easy to see (and
proved in Casselman and Shalika~\cite{CS}) that
$$T_1F_k=q^{-k/2}\,\gamma_1^k,\qquad
T_{w_0}F_k=(1-q^{-1}\gamma_1\gamma_2^{-1})
(1-\gamma_1\gamma_2^{-1})^{-1}\,q^{-k/2}\,\gamma_2^k. \tag2.5$$
Thus
$$W(F_k)=W(f_1)\,q^{-k/2}\,\gamma_1^k+
(1-q^{-1}\gamma_1\gamma_2^{-1})
(1-\gamma_1\gamma_2^{-1})^{-1}\,q^{-k/2}\,\gamma_2^k\,W(f_{w_0}).$$
Now using (2.3), we find that
$$\multline(1-q^{-1})W(\zeta_n)=(1-q^{-1/2}\tau\gamma_1^{-1})
W(f_1)\,q^{-k/2}\,\gamma_1^k+\\(1-q^{-1/2}\tau\gamma_2^{-1})
(1-q^{-1}\gamma_1\gamma_2^{-1})
(1-\gamma_1\gamma_2^{-1})^{-1}\,q^{-k/2}\,\gamma_2^k\,W(f_{w_0}).
\endmultline$$

We may compute $W(f_{w_0})$ explicitly. By definition this equals
$$\int_{F^\times}f_{w_0}\!\left(\pmatrix &1\\1\endpmatrix
\pmatrix1&1\\&1\endpmatrix\pmatrix a\\&1\endpmatrix\right)
\sigma(a)^{-1}\,d^\times a,$$
and it follows from (2.4) that the integrand here equals
$$\phi_\xi\!\left(\pmatrix &1\\1\endpmatrix
\pmatrix1&1\\&1\endpmatrix\pmatrix a\\&1\endpmatrix\right)$$
if $|a|\ge1$, zero otherwise. Thus
$$W(f_{w_0})=\int_{|a|\ge1}\phi_\xi\!\left(\pmatrix &1\\1\endpmatrix
\pmatrix1&1\\&1\endpmatrix\pmatrix a\\&1\endpmatrix\right)
\sigma(a)^{-1}\,d^\times a.$$
If $|a|=q^{k}$, $k\ge0$, then the integrand is readily evaluated using
(2.1), and equals
$\tau^k\,\gamma_2^{-k}\,q^{-k/2}$, so, assuming (1.6), the last integral is
absolutely convergent and equals $(1-q^{-1/2}\tau\gamma_2^{-1})^{-1}$.
Thus
$$\multline(1-q^{-1})W(\zeta_n)=(1-q^{-1/2}\tau\gamma_1^{-1})
W(f_1)\,q^{-k/2}\,\gamma_1^k+\\
(1-q^{-1}\gamma_1\gamma_2^{-1})
(1-\gamma_1\gamma_2^{-1})^{-1}\,q^{-k/2}\,\gamma_2^k.
\endmultline$$

We have not yet proved the Theorem since we have not evaluated $W(f_1)$.
Note however that the Theorem will follow if we prove that $\Wa_\xi$ is
invariant under the interchange of $\gamma_1$ and $\gamma_2$; indeed, we
may assume without loss of generality that $\xi$ is regular, so
$\gamma_1\ne\gamma_2$. Then, because the two functions $q^{-k/2}\gamma_1^k$
and $q^{-k/2}\gamma_2^k$ of $k$ are linearly independent, the unknown value
of $W(f_1)$ will be determined.

We will show that $\Wa_\xi(1)=1$. This is sufficient: indeed, recalling
Waldspurger's theorem on the uniqueness of the model, there is a unique
spherical vector in this unique model which is normalized to equal $1$ at
$g=1$, if $\Wa_\xi(1)=1$ then clearly $\Wa_\xi$ must be this vector; then,
since the isomorphism class of the representation $\pi_2=\text{Ind}(\xi)$
is unchanged when we interchange $\gamma_1$ and $\gamma_2$, $\Wa_\xi$ is
thus invariant under this interchange. By definition
$$W\!a_\xi(1)=\int_{F^\times}\phi_\xi\!\left(\pmatrix &1\\1\endpmatrix
\pmatrix1&1\\&1\endpmatrix\pmatrix a\\&1\endpmatrix\right)
\sigma(a)^{-1}\,d^\times a,$$
and in this integral, the integrand is constant when $a$ has constant
valuation. (This is not true for $W\!a_\xi(\eta_k)$, which is the reason
for the somewhat elaborate proof of this Theorem!) Once again
applying (2.1), we find that the
integrand equals $\tau^{-k}\gamma_1^{k}q^{-k/2}$ if $|a|=q^{-k}$, $k\ge0$, or
$\tau^{-k}\gamma_2^kq^{k/2}$ if $k<0$.  Assuming (1.6) holds,
it is then simple to sum the two
geometric series and check that $\Wa_\xi(1)=1$.

This completes the proof of Theorem 1.1\smudge

We turn now to the proof of Theorem~1.2. Fix an element
$f\in C^\infty\big((B_2\cap K_2)\backslash K_2\big)$, and extend
$f$ to $\Ind(\xi)$.  Then, since $f$ is
locally constant, it follows from (2.1) that the
integrand on the right in (1.5) equals
$$\cases |a|^{1/2}\xi_1(a)\sigma^{-1}(a)\,f\!\pmatrix -1&0\\0&1\endpmatrix&
\text{if $|a|$ is sufficiently small;}\\
|a|^{-1/2}\xi_2(a)\sigma^{-1}(a)\,f\!\pmatrix0&1\\1&0\endpmatrix&
\text{if $|a|$ is sufficiently large.}\endcases$$
Hence the integral (1.5) is equal to an integral over a compact set
plus two integrals giving geometric series, whose values have analytic
continuation to the region $\Lambda$. This gives the analytic continuation
of $W$, and it only remains to be seen that it represents a Waldspurger
functional. Thus we must show that with $t\in T_2^{\text s}(F)$,
$W(\pi_2(t)\,f-\sigma(t)\,f)=0$. It is clear that this is true when
(1.6) is satisfied, and that the left side is analytic, so this is
true for all $(\gamma_1,\gamma_2,\tau)\in\Lambda$. This completes the
proof of Theorem~1.2. \smudge

\medbreak
\subheading{3. Proof of the Analytic Continuation and Functional Equation
for the Bessel Model}
In this Section we shall prove (most of) Theorems 1.3 and 1.4.  The
proof of the continuation in the parameters $\alpha_i$ and of the
functional equation is based on homomorphisms from $GL(2,F)$ into $G$,
similarly to Jacquet's proof of the analytic continuation and
functional equation of the Whittaker functions on Chevalley groups
\cite{Ja}. Jacquet's method suffices to give most of the functional
equations, but an extra step is needed (different in the nonsplit and
split cases).  One then applies Hartog's theorem. At this point, we
will have proved Theorems~1.3 and 1.4 except for one point, namely the
meromorphic continuation in the split case outside the region
$q^{-1/2}<\min(|\beta|,|\beta|^{-1})$. This
meromorphic continuation follows from the explicit formula in
Theorem~1.6, or from Theorem~1.7; so for this minor point, the proof
will be completed in subsequent sections.

Our proof will follow to the extent possible the notation and
organization of \cite{BFG}, where a similar method was used to study
another unique functional.

First, we recall two well-known lemmas.  To give the first, let
$\omega\in\Omega$ be a Weyl group element represented by the permutation
matrix $w$ (we shall frequently abuse the notation and write
$w\in\Omega$). Let $N$ denote the full subgroup of upper triangular
unipotent matrices in $G$.  Given a character $\chi$ as above, let
$^w\!\chi$ be the character satisfying $^w\!\chi(a)=\chi(w^{-1}aw)$ for all
diagonal $a\in G$.  Define the intertwining operator
$T_w\colon\Ind(\chi)\to\Ind(^w\!\chi)$ by the integral
$$(T_w\Psi)(g)=\il{N_w\bs N}\Psi(w^{-1}ng)\,dn,\tag3.1$$
where $N_w=N\cap wNw^{-1}$.  Given an unramified character $\chi$ as
above, let $a_\alpha$ be the diagonal matrix in $G$
$$a_\alpha=(\alpha_1,\cdots,\alpha_n,1,\alpha_n^{-1},\cdots,
\alpha_1^{-1}).$$
Order the roots of $SO(2n+1,\C)$ so that $N$ corresponds to the positive
ones.  Then one has (\cite{Ca2}, Section 6.4)

\proclaim{Lemma~3.1}  The intertwining integral (3.1) is absolutely
convergent if $|r(a_\alpha)|<1$ for all positive roots $r$ of $SO(2n+1,\C)$
such that $w(r)<0$.  Moreover, $T_w$ varies holomorphically with $\chi$
and has meromorphic continuation to the space of all unramified characters.
\endproclaim

The second lemma concerns the $G_2=GL(2,F)$ Whittaker function.
As in Section 1,
given two unramified quasicharacters $\xi_1$, $\xi_2$ of $F^\times$,
$\xi_i(\varpi)=\gamma_i$,
define a character $\xi$ of the standard Borel subgroup of $G_2$
by  equation (1.3).
Let $\pi_2=\Ind(\xi)$ be the normalized induced
representation of $G_2$,
and $\phi_\xi$ be the $K_2$-fixed vector such that $\phi_\xi(I_2)=1$.
If $|\gamma_1|>|\gamma_2|$,
let $W\!h_\xi$ be the Whittaker function
$$W\!h_\xi(g)=\il{F}\phi_\xi\left(\pmatrix&1\\1&\endpmatrix
\pmatrix 1&x\\&1\endpmatrix g\right)\,\psi(x)\,dx.\tag3.2$$
Then one has
\proclaim{Lemma 3.2}
\roster
\item The Whittaker function $W\!h_\xi$, originally
defined by the integral (3.2) when $|\gamma_1|>|\gamma_2|$, has
a meromorphic continuation to all nonzero complex $\gamma_1,\gamma_2$.
Moreover, the function
$$\Wh_\xi(g)=(1-\gamma_1\gamma_2^{-1}q^{-1})^{-1}W\!h_\xi(g)$$
is holomorphic in $(\C^\times)^2$ and is invariant under the
interchange of $\gamma_1$ and $\gamma_2$.
\item Let $\gamma_1'=(1+\epsilon)\,\max(|\gamma_1|,|\gamma_2|)$,
$\gamma_2'=(1+\epsilon)^{-1}\,\min(|\gamma_1|,|\gamma_2|)$,
where $\epsilon\geq0$ is chosen so that  $\gamma_1'\neq\gamma_2'$.
Define corresponding unramified quasicharacters
$\xi_i'$ for $i=1,2$ by $\xi_i'(\varpi)=\gamma_i'$.
Then
$$|W\!h_\xi(g)|\ll\left|\phi_{\xi'}(g)\right|$$
uniformly in $g$,
where $\phi_{\xi'}$ is the normalized $K_2$-fixed
vector in $\Ind(\xi_1',\xi_2')$.
\endroster
\endproclaim

Just as the continuation of the Waldspurger function $W\!a_\xi$
is deduced from its evaluation,
Lemma 3.2 may be deduced from the explicit evaluation of the Whittaker
function $W\!h_\xi$.  Assuming $|\gamma_1|>|\gamma_2|$,
a computation shows that the integral (3.2) is zero if
$g=\left(\smallmatrix \varpi^k&0\\0&1\endsmallmatrix\right)$
with $k<0$, and is given by
$$W\!h_\xi\left(\pmatrix \varpi^k&0\\0&1\endpmatrix\right)=
q^{-k/2}(1-\gamma_1\gamma_2^{-1}q^{-1})\,
\frac{\gamma_1^{k+1}-\gamma_2^{k+1}}{\gamma_1-\gamma_2}$$
if $k\geq0$.  The Lemma then follows.

We give next a third Lemma, concerning the convergence of the
integrals which arise in the consideration of the split case.
Let $\chi$ be an unramified character
as above.  Let $w\in\Omega$ be a given permutation, and
factor $w^{-1}$ as $w^{-1}={w'}^{-1}(n,n+2)$.
Set $U_w=U\cap w^{-1}Uw$.  Also, define the integer $g$
and the complex number $\alpha$ by the equations
$$\delta_B\left({w'}^{-1}t(\varpi)w'\right)=q^{-g},\qquad
\chi\left({w'}^{-1}t(\varpi)w'\right)=\alpha.$$

\proclaim{Lemma~3.3}  For $\Psi\in\Ind(\chi)$, the integral
$$\il{F^\times}\il{U_w\bs U}\Psi(w^{-1}n(1)\,u\,t(a))\,
\lambda^{-1}(a)\,du\,d^\times a\tag3.3$$
is absolutely convergent provided
$$|\alpha|<q^{g/2}\,\min(|\beta|,|\beta^{-1}|)$$
and provided that the intertwining integral $T_{w'}(\Psi)$
converges absolutely.
\endproclaim

\demo{Proof} We may interchange $u$ and
$t(a)$ in the above integral. We have
$w^{-1}n(1)\,t(a)\,u={w'}^{-1}(n,n+2)\,n(1)\,t(a)\,u$.
A computation similar to (2.1) shows that
$$(n,n+2)\,n(1)\,t(a)=\cases t(a)\,n(a^{-1})\,\kappa_1(a)&\text{if
$|a|_F\leq1$}\\t(a^{-1})\,\kappa_2(a)&\text{if $|a|_F>1$},
\endcases$$
with $\kappa_1(a)$, $\kappa_2(a)$ in $K$ of the form
$$\pmatrix I_{n-1}&&\\&*&\\&&I_{n-1}\endpmatrix.$$
Since matrices of this form normalize $U$, we may move
the $\kappa_i(a)$ to the right in the integral (3.3).
Using the condition
$\Psi\in\Ind(\chi)$ and comparing with the definition
(3.1) of the intertwining operator $T_{w'}$, one sees that the
integral (3.3) is absolutely bounded by the sum of the
two integrals
$$\il{|a|_F\leq1}|\alpha|^{\ord(a)}\,|\lambda^{-1}(a)|\,
|a|_F^{g/2}\,|(T_{w'}\Psi)(\kappa_1(a))|\,d^\times a$$
and
$$\il{|a|_F>1}|\alpha|^{-\ord(a)}\,|\lambda^{-1}(a)|\,
|a|_F^{-g/2}\,|(T_{w'}\Psi)(\kappa_2(a))|\,d^\times a.$$
Since $(T_{w'}\Psi)(g)$ is a locally constant function, it is absolutely
bounded on $K$.  Hence both integrals are bounded by
geometric series, which converge under the hypotheses
of the Lemma. \smudge
\enddemo

This proof may also be rephrased by using the isomorphism of
$PGL(2,F)$ with $SO(3,F)$ to write the integral (3.3) as the
Waldspurger integral of an intertwining integral, and then applying
Lemmas 2.1 and Lemma 3.1.

Observe that applying Lemma 3.3 with $w=w_0$, so that $w'=w_1$, $g=1$
and $\alpha=\alpha_n$, one finds that the integral (1.14) is
absolutely convergent in the region (1.15), as claimed.

To prove Theorems~1.3 and 1.4 we now proceed in two stages.  First,
using the standard functional equation for the $GL(2)$~$p$-adic
Whittaker function (Lemma 3.2), we shall establish functional
equations under the transpositions $(j,j+1)\in\Omega$, $1\leq j\leq
n-1$.  The proofs of these results also give the analytic continuation
to certain unions of Weyl chambers properly larger than the original
region of convergence.  Then, using Theorem 1.1 in the split case and
the invariance of the $G_2$-spherical function in $\Ind(\xi_1,\xi_2)$
under the interchange of $\xi_1$ and $\xi_2$ in the
nonsplit case, we obtain a functional equation for the Weyl group
element taking $\alpha_{n}$ to $\alpha_{n}^{-1}$ and fixing the rest.
Since these two elements generate $\Omega$, these steps imply that
$\H_{\chi}$ is invariant under $\Omega$.

To obtain the functional equations under the transpositions $(j,j+1)$,
$1\leq j\leq n-1$, let $\iota_j$ be the
embedding of $GL(2,F)$ into $G$ given by
$$g\mapsto\pmatrix I_{j-1}&&&&\\&\det(g)^{-1}g&&&\\&&I_{2n-2j-1}&&\\
&&&g^\sharp&\\&&&&I_{j-1}\endpmatrix,$$
where $g^\sharp=\left(\smallmatrix -1&\\&1\endsmallmatrix\right)g
\left(\smallmatrix -1&\\&1\endsmallmatrix\right)$.
Let $w$ denote $w_0$ if $T$ is split, and $w_1$ if $T$ is nonsplit.
Let $v_1=(j,j+1)(2n-j+1,2n-j+2)$, and factor $w=v_1v_2$.
Let $B_2$ and $U_2$ be, as above, the subgroups of $GL(2,F)$
consisting respectively of upper triangular and of upper triangular
unipotent matrices,
and let $B_j$ and $U_j$ be the subgroups of $G$ given by
$$B_j=\iota_j(B_2),\qquad
U_j=v_2^{-1}\,\iota_j(U_2)\,v_2.$$
For $u_j\in U_j$, write
$v_2\,u_j=\iota_j(\hat{u}_j)\,v_2$.
Let $U'_j$ and $B'_j$ be the complementary
subgroups in $U$ and $B$ to $U_j$ and $B_j$, respectively,
so that $U=U_jU'_j$ and $B=B_jB'_j$ (uniquely).
Note that $v_1$ and $\iota_j(U_2)$ normalize $B_j'$
and fix $\Phi_\chi|_{B_j'}$.
Let $u\in U$, $t\in T$, $a\in G$.  Factor $u=u_ju'_j$,
$u_j\in U_j$, $u'_j\in U'_j$.  Applying the Iwasawa decomposition,
for $g\in G$, we may write $v_2u'_jtg=b'_jb_j\kappa$, with
$b_j\in B_j$, $b'_j\in B'_j$, $\kappa\in K$ (we suppress the dependence
of $b'_j$ and $b_j$ on $u'_j$, $t$, and $g$ from the notation).

Suppose first that $T$ is nonsplit. Then with the above notation we have
$$\align
w_1utg&=v_1v_2u_ju'_jtg\\
&=v_1\iota_j(\hat{u}_j)\,v_2 u_j'tg\\
&=v_1\iota_j(\hat{u}_j)\,b'_jb_j\kappa.
\endalign$$
Since $\Phi_\chi$ is right $K$-invariant, we may thus
express $H_\chi(g)$ as the iterated integral
$$H_\chi(g)=\il{T(\O)}\il{U'_j}\left[\;\il{U_j}
\Phi_\chi(v_1\iota_j(\hat{u}_j)b_j)\,
\theta_S(u_j)^{-1}\,du_j\right]
\Phi_\chi(b'_j)\,\theta_S(u'_j)^{-1}\,du'_j\,dt.\tag3.4$$
However for $h\in G$, it follows from (1.11) that
the function $\phi\colon GL(2,F)\to\C$ given by
$\phi(a)=\Phi_\chi(\iota_j(a)h)$ is in
the space $\Ind(\chi_{j+1}^{-1}\mu^{j-n},
\chi_j^{-1}\mu^{j-n})$, where
$\mu$ is the quasicharacter of $F^\times$ given by $\mu(x)=|x|_F$.
Accordingly the inner integral in (3.4) is a constant
multiple of the $GL(2)$ Whittaker function
associated to this representation.
Applying Lemma 3.2, we obtain the analytic continuation of the function
$${1\over 1-\alpha_j\alpha_{j+1}^{-1}q^{-1}}H_\chi(g)$$
to the region $C_j$ of $(\C^\times)^{n}$ defined by
the inequalities
$$\multline|\alpha_1|<\cdots<|\alpha_{j-1}|<\min(|\alpha_j|,|\alpha_{j+1}|)
\leq\max(|\alpha_j|,|\alpha_{j+1}|)\\<|\alpha_{j+2}|<\cdots<|\alpha_{n-1}|
<\min(|\alpha_{n}|,|\alpha_{n}^{-1}|)\endmultline\tag3.5$$
if $1\leq j<n-1$, and defined by the inequalities
$$|\alpha_1|<|\alpha_2|<\cdots<|\alpha_{n-2}|< \min(|\alpha_{n-1}|,
|\alpha_{n}|)\leq \max(|\alpha_{n-1}|,|\alpha_{n}|)<1\tag3.6$$
if $j=n-1$, and the invariance of this function under $(j,j+1)$ there.
The region $C_j$ is obtained by replacing the inner
integral by the estimate given in Lemma 3.2, part (2),
and comparing with the intertwining operator $T_{v_2^{-1}}$,
whose convergence is given in Lemma 3.1.  The difference
in inequalities is due to this comparison; note that
for $1\leq j<n-1$, $v_2$ is a product of transpositions
$$v_2=(1,2n+1)\cdots(j-1,2n-j+3)(j,2n-j+1)(j+1,2n-j+2)
(j+2,2n-j)\cdots(n-1,n+3),$$
while for $j=n-1$, $v_2$ includes a four-cycle
$$v_2=(1,2n+1)\cdots(n-2,n+4)(n-1,n+2,n+3,n).$$
Also, let us remark that if $1\leq j<n-1$, then the region $C_j$
properly contains the original region of convergence $C_0$,
while if $j=n-1$, Hartog's theorem gives at once the
analytic continuation of $H_\chi$ to $C_0\cup C_{n-1}$.

Suppose now that $T$ is split.  Let $u$, $u_j$, $u'_j$,
$\hat{u}_j$ be as above; note that $U_j=\iota_j(U_2)$.  We have
$$\align
w_0n(1)ut(a)&=v_1v_2u_j (u_j^{-1}n(1)u_j)u'_jt(a)\\
&=v_1\iota_j(\hat{u}_j)v_2(u_j^{-1}n(1)u_j)u'_jt(a).
\endalign$$
Moreover, for
$u_j=u_j(x):=\iota_j\left(\left(\smallmatrix 1&x\\&1
\endsmallmatrix\right)\right)$, a calculation shows that
$u_j^{-1}n(1)u_j=n(1)u_j''$ with $u_j''\in U'_j$ and
$\theta_S(u_j)^{-1}\theta_S({u_j''})=\psi(-x)$.
Using the Iwasawa decomposition, write
$v_2n(1)u'_jt(a)g=b'_jb_j\kappa$, with
$b_j\in B_j$, $b'_j\in B'_j$, $\kappa\in K$
(once again we suppress from the notation the dependence
of $b'_j$ and $b_j$ on $u'_j$, $t(a)$, and $g$).
Then the integral (1.14) representing $H_\chi(g)$ becomes
$$H_\chi(g)=\il{F^\times}\il{U'}\left[\;\il{F}
\Phi_\chi(v_1 u_j(-x) b_j)\,\psi(-x)\,dx\right]
\Phi_\chi(b'_j)\,\theta_S(u'_j)^{-1}
\,\lambda^{-1}(a)\,du\,d^\times a.$$
The inner integral is once again a $GL(2)$ Whittaker function.
Applying Lemma 3.2, we obtain the analytic continuation
of the function
$${1\over 1-\alpha_j\alpha_{j+1}^{-1}q^{-1}}H_\chi(g),$$
to the region $C'_j$ of $(\C^\times)^{n}$ defined by
requiring the inequalities (3.5) and in addition the inequality
$$|\alpha_n|<q^{1/2}\min(|\beta|,|\beta^{-1}|)$$
if $1\leq j<n-1$, and defined by requiring the inequalities
(3.6) and in addition the inequality
$$\min(|\alpha_{n-1}|,|\alpha_n|)<q^{3/2}\min(|\beta|,|\beta^{-1}|)$$
if $j=n-1$.
We also obtain the invariance of this function under $(j,j+1)$ there.
The region $C'_j$ is obtained by replacing the inner
integral by the estimate given in Lemma 3.2, part (2),
and applying Lemma 3.3.  For later use, we denote the original
region of convergence given by (1.15) as $C_0'$.

It remains to obtain the functional equation under the
interchange $\alpha_{n}\leftrightarrow\alpha_{n}^{-1}$.
To obtain this, let $\iota_n$ denote the homomorphism of
$GL(2,F)$ into $G$ given by
$$\multline\iota_n\left(\pmatrix a&b\\c&d\endpmatrix\right)
=\\{1\over ad-bc}
\pmatrix (ad-bc)I_{n-1}&&&&\\&a^2&ab&-b^2/2&\\&2ac&ad+bc&-bd&\\
&-2c^2&-2cd&d^2&\\&&&&(ad-bc)I_{n-1}\endpmatrix.\endmultline
\tag3.7
$$

Consider first the case that $T$ is nonsplit.  Since $w_1$ fixes
$T$, the integral (1.12) representing $H_\chi(g)$
is expressed as an iterated integral
$$H_\chi(g)=\il{U}\left[\il{T(\O)}\Phi_\chi(tw_1ug)\,dt\right]\,
\theta_S(u)^{-1}\,du.\tag3.8$$
Now for $h\in G$, the function $\phi\colon GL(2,F)\to\C$ given by
$\phi(a)=\Phi_\chi(\iota_n(a)h)$ is in the space $\Ind(\chi_n,\chi_n^{-1})$.
Recalling that $\iota_n$ gives the isomorphism
between $PGL_2(F)$ and $SO(3,F)$, one sees that the inner
integral in (3.8) is a multiple of the nonsplit Waldspurger
functional on $PGL_2(F)$, evaluated on a suitable right-translate
of the spherical vector in this induced space.
As explained in Section 1, it is thus
a multiple of the spherical function in $\Ind(\chi_n,\chi_n^{-1})$.
By the invariance of this spherical function under the interchange
of $\chi_n$ and $\chi_n^{-1}$, we conclude
that $H_\chi(g)$ is invariant under
$\alpha_n\leftrightarrow\alpha_n^{-1}$ in the domain (1.13).

Consider instead the case that $T$ is split.
Factor $w_0=w_1w_2$ with $w_2=(n,n+2)$.
In this case the integral (1.14) is an iterated integral of the form
$$H_\chi(g)=\il{U}\left[\,\,\il{F^\times}\Phi_\chi(w_2\,n(1)\,t(a)\,w_1\,u\,g)
\,\lambda^{-1}(a)\,d^\times a\right]\,\theta_S(u)^{-1}\,du.$$
Once again using the homomorphism $\iota_n$, one sees that
the inner integral is a multiple of the split Waldspurger functional
computed in Section 2, applied to a function in
$\Ind(\chi_n,\chi_n^{-1})$.
Note that this integral is absolutely convergent
provided $|\alpha_n|<q^{1/2}\min(|\beta|,|\beta^{-1}|)$.
Applying Theorem 1.1 and arguing as above,
we obtain the invariance of the function
$${{(1-\alpha_n\beta q^{-1/2})(1-\alpha_n\beta^{-1} q^{-1/2})}
\over{(1-\alpha_n^2 q^{-1})}}
\,H_{\chi}(g),$$
under $\alpha_n\leftrightarrow\alpha_n^{-1}$ in the domain (1.15).

To complete the proof of Theorem~1.3,  let $C$ be the set of
$\alpha=(\alpha_1,\cdots,\alpha_{n})\in
(\C^\times)^{n}$ such that at most one of the equalities
$|\alpha_i|=|\alpha_j|$, $1\leq i < j \leq n$, and
$|\alpha_i|=|\alpha_j|^{-1}$, $1\leq i \leq j \leq n$,
is satisfied.
It may be deduced from Hartog's theorem that any analytic function
on $C$ can be extended to an analytic function on $(\C^\times)^{n}$, and so
it is sufficient to extend the function $\H_{\chi}(a)$ to the domain $C$ in
such a way that the corresponding functional equations are satisfied on
$C$. Now if $\alpha\in C$, then there is an $\omega\in\Omega$ such that
$\omega\alpha$ is in $C_k$ for some $k$, $0\leq k\leq n-1$.
We then define $\H_{\chi}(a)$ to equal $\H_{^\omega\!\chi}(a)$.  This is
well defined by the above discussion. It is apparent that this extends
$\H_{\chi}(a)$ to an analytic function of $C$ satisfying the corresponding
functional equations, as required. This completes the proof of
Theorem~1.3. \smudge

The proof of Theorem~1.4 is similar.  Suppose first that
$\lambda$ is chosen so that
$q^{-1/2}<\min(|\beta|,|\beta|^{-1})$.  Then using the regions $C'_j$
in place of $C_j$, the argument in
the paragraph above gives the analytic continuation of $\H_{\chi}$
to $\alpha\in(\C^\times)^n$, and the functional equation there.
This completes the proof of Theorem~1.4,
except, as we have noted, for the point of the meromorphic
continuation to all $\alpha_i$ and $\beta$, and this point is a
consequence of the explicit formula in Theorem~1.6 or of Theorem~1.7.
\smudge

\goodbreak
\subheading{4.  Explicit Formulas for the Bessel Model}
The evaluation of the Bessel model given in Theorems 1.5 and
1.6 is obtained by applying the method of Casselman and Shalika
\cite{Ca1}, \cite{CS}.  This method was also used in
Section 2.
A similar evaluation is carried out \cite{BFG}; however,
the particulars are different in the case at hand,
and the cases $T$ nonsplit and $T$ split are
once again different from each other.
Throughout the proofs we shall assume that
$\chi$ is regular, i.e. $^\omega\chi\neq\chi$ for all nonidentity
$\omega\in \Omega$. The general case then
follows from the analytic continuation given in Theorems~1.3 and 1.4.

Let $\B$ be the Iwahori subgroup of
$K$, consisting of integral matrices which are upper triangular
invertible $\bmod \varpi\O$. It follows from the Iwasawa decomposition
of $G$ and the Bruhat decomposition over $\O/\varpi\O$
that the space of right-Iwahori-fixed vectors
$\Ind(\chi)^{\B}$ is $|\Omega|$-dimensional; moreover, a basis is
given by the functions $\phi_w$ defined by
$$\phi_w(bw^{\prime\,-1}b_1)=\cases \Phi_\chi(b)&
\text{if $w=w\dash$;}\\0&\text{otherwise,}\endcases$$
where $b\in B$, $w\in \Omega$, $b_1\in\B$
(note that this differs from Casselman's notation in
\cite{Ca1}, but is consistent with \cite{BFG}).

If $\chi$ is regular, then it is shown in \cite{Ca1},
Section 3, that the
linear functionals on $\Ind(\chi)^{\B}$ given by
$f\mapsto (T_wf)(I_{2n+1})$ are linearly independent.
Here $T_w$ is the intertwining operator defined in (3.1).
Let $f_w$, $w\in\Omega$, be the dual basis, characterized by
$$(T_wf_{w\dash})(I_{2n+1})=\cases 1&\text{if $w=w\dash$;}
\\0&\text{otherwise.}\endcases$$

Suppose first that $T$ is nonsplit.
Let $d=d_k$ be as given in Section 1, and suppose that all $k_i\geq0$.
Let $F_d$ be the function
$$F_d(g)=\il{T(\O)}\il{U\cap K}\Phi_\chi(gutd)\,du\,dt.$$
Clearly $F_d\in\Ind(\chi)$. Let us show that $F_d$ is
right Iwahori invariant; then we may write
$$F_d(g)=\sum_{w\in \Omega}R(d,w;\chi)\,f_w.\tag4.1$$

To prove the Iwahori invariance, we first note that by Jacquet's first lemma
(see the proof of Proposition~2.5 of Casselman~\cite{Ca1}),
$F_d$ is (right) invariant under the group $\Sigma$ of elements of
$SO(2n+1,\O)$ whose lower entries are divisible by $\varpi$, and whose
middle $3\times 3$ block is the identity; and also, if $\iota_n$ is as
in (3.7), it is clear that $F_d$ is invariant by
$\iota_n\big(T_2^{\text a}(\O)\big)$. Thus
$$F_d(bg\sigma)=(\delta_B^{1/2}\chi)(b)F_d(g)\tag4.2$$
when $b\in B(F)$ and $\sigma\in\iota_n\big(T_2^{\text a}(\O)\big)\Sigma$.
We will deduce invariance by the group $\iota_n(K_2)\Sigma$ (which contains
the Iwahori subgroup) from the fact that the canonical map
$$B\backslash G/\iota_n(T_2^{\text a}(\O))\Sigma\to
B\backslash G/\iota_n(K_2)\Sigma\tag4.3$$
is a bijection. To see this, we note that since $\iota_n(K_2)\Sigma$ contains
the Iwahori subgroup, every double coset in
$B\backslash G/\iota_n(K_2)\Sigma$ contains a Weyl group element
$w$; so what we must show is that
$Bw\iota_n(T_2^{\text a}(\O))\Sigma=Bw\iota_n(K_2)\Sigma$.
We consider an element $bw\iota_n(k)\sigma$ of the right side,
where $b\in B$, $k\in K_2$ and $\sigma\in\Sigma$. We can write
$k=\beta^+t^+=\beta^-t^-$ where $t^\pm$ are in $T_2^{\text a}(\Cal O)$,
$\beta^+$ is upper triangular and $\beta^-$ is lower triangular. Then
one of $w\iota_n(\beta^\pm) w^{-1}$ is upper triangular, so
$bw\iota_n(k)\sigma=b'w\iota_n(t^\pm)\sigma$ where
$b'=bw\iota_n(\beta^\pm)w^{-1}$,
which shows that this element lies in $Bw\iota_n(T_2^{\text a}(\O))\Sigma$.
Thus (4.3) is a bijection. We may now prove the Iwahori invariance of
$F_d$. If $g\in G$ and $\sigma\in\iota_n(K_2)\Sigma$, we write
$g\sigma=bg\sigma'$ with $b\in B$ and
$\sigma'\in \iota_n\big(T_2^{\text a}(\O)\big)\Sigma$. We note that $b$
is conjugate to $\sigma{\sigma'}^{-1}\in K$, and so the eigenvalues of $b$
are units, and $\chi(b)=1$. It thus follows from (4.2) that
$F_d(g\sigma)=F_d(bg\sigma')=F_d(g)$.

Let us compute the coefficient $R(d,w;\chi)$. It equals
$$\align
(T_wF_d)(I_{2n+1})&=\il{N_w\bs N}\il{T(\O)}\il{U\cap K}
\Phi_\chi(w^{-1}nutd) \,du\,dt\,dn\\
&=\il{T(\O)} T_w(\Phi_{\chi})(td)\,dt.
\endalign$$
However, let $\Phi^+$ denote the set of positive roots of
$SO(2n+1,\C)$, and if $r\in\Phi^+$, let $\iota_r$ be the
corresponding embedding of $SL(2,F)$ into $G$, and
$$a_r=\iota_r\pmatrix \varpi&\\&\varpi^{-1}\\\endpmatrix.$$
Then it is shown in \cite{Ca1}, Theorem 3.1, that
$$T_w(\Phi_{\chi})=c_w(\chi)\,\Phi_{^w\!\chi},$$
where $\Phi_{^w\!\chi}$ is the standard nonramified vector
in $\Ind(^w\!\chi)$, and
the coefficient $c_w(\chi)$ is given by
$$c_w(\chi)=\prod_{\scriptstyle r\in\Phi^+\atop\scriptstyle w(r)<0}
\left({1-q^{-1}\,\chi(a_r)\over1-\chi(a_r)}\right).\tag4.4$$
We find that
$$R(d,w;\chi)=c_w(\chi)\,\sigma_{^w\!\chi}(d),$$
where
$$\sigma_{^w\!\chi}(d)=\il{T(\O)}
\Phi_{^w\!\chi}(td)\,dt.$$

Let $\iota_n$ be the homomorphism from $GL(2,F)$ into $G$ given by
equation (3.7).
Factor $d=d'\iota_n(d'')$ (all these matrices depending on the $k_i$) with
$$d'=\text{diag}(\varpi^{k'_1},\cdots,\varpi^{k'_2},1,1,1,
\varpi^{-k'_2},\cdots,\varpi^{-k'_n})$$
and
$$d''=\pmatrix \varpi^{k_1}\\&1\endpmatrix.$$
Then one sees that
$$\sigma_{^w\!\chi}(d)=
{}^w\!\chi\delta_B^{1/2}(d')\!\il{T(\O)}
\Phi_{^w\!\chi}(t\,\iota_n(d''))\,dt.$$
Now if
$$^w(\chi_1,\cdots,\chi_{n})=
(\chi\dash_1,\cdots,\chi\dash_{n}),$$
then the function $g\mapsto\Phi_{^w\!\chi}(\iota_n(g))$
lies in the induced space
$\Ind(\chi_{n}\dash,{\chi_{n}\dash}^{-1})$.
Recalling that $\iota_n$ identifies $PGL(2,F)$ with $SO(3,F)$
and applying the analysis of the nonsplit Waldspurger
functional in Section 1, one concludes that
this last integral is equal to the spherical function
on $GL(2,F)$ in this induced space, evaluated at $d''$.
This value is given by the Macdonald
formula (see \cite{Ca1}).  Substituting this formula, we obtain
$$\multline\sigma_{^w\!\chi}(d)=
(1+q^{-1})^{-1}\,q^{e_k}\,
\prod_{i=2}^{n}\chi\dash_{n+1-i}(\varpi^{k'_i})\\
\times{{(\chi_n\dash(\varpi)-{\chi_n\dash}^{-1}(\varpi)q^{-1})\,
\chi_n\dash(\varpi)^{k_1}-
({\chi_n\dash}^{-1}(\varpi)-{\chi_n\dash}(\varpi)q^{-1})\,
\chi_n\dash(\varpi)^{-k_1}
}\over{\chi_n\dash(\varpi)-{\chi_n\dash}^{-1}(\varpi)}}.
\endmultline\tag4.5$$
This completes the evaluation of $R(d,w;\chi)$.

If $f\in\Ind(\chi)$ with $\chi$ dominant, let
$$L(f)=\il{U}f(w_1u)\,\theta_S(u)^{-1}\,du.$$
This integral converges absolutely by comparison with $T_{w_1}(f)$.
Then, by (1.12) and (4.1), we have
$$H_\chi(d)=L(F_d)=\sum_{w\in \Omega}R(d,w;\chi)\,L(f_w).\tag4.6$$
It is clear from this formula and our evaluation of $R(d,w;\chi)$ that
$h(k_1,\cdots,k_{n})$ is a linear combination of functions of the
$\alpha_i$ which lie in exactly one $\Omega$-orbit.
To complete the proof of Theorem~1.5, we
shall compute the full contribution to $h(k_1,\cdots,k_n)$ of
the terms on the right of (4.6) with $w=w_0$, $w=w_1$.  By
the computation of $R(d,w;\chi)$ given above, no other $w\in\Omega$
contributes a rational function of the form
$$\alpha_1^{-k_n'}\cdots\alpha_n^{-k_1'}$$
times a function independent of the $k_i$.  Hence, by Theorem 1.3, we
will have computed one piece of the $\Omega$-orbit, and the final
formula is then the $\Omega$-symmetric sum of these terms, taking
into account the normalizing factor in the functional equation of
Theorem 1.3.

To compute the contribution of the $w=w_0,w_1$ terms to (4.6), we need the
following lemma relating the two bases for the Iwahori fixed vectors,
$\{\phi_w\}$ and $\{f_w\}$.

\proclaim{Lemma~4.1}
\roster
\item $\phi_{w_0}=f_{w_0}$.
\item $\phi_{w_1}=f_{w_1}+\frac{(1-q^{-1})\alpha_n^2}{1-\alpha_n^2}f_{w_0}$.
\item $L(f_{w_0})=0$.
\item $L(f_{w_1})=1$.
\endroster
\endproclaim

\demo{Proof} For $w'\in\Omega$ we may write
$$\phi_{w'}(g)=\sum_{w\in\Omega}c_\chi(w,w')\,f_w(g).\tag4.7$$
where the coefficients $c_\chi$ are given by
$$\align
c_\chi(w,w')&=T_w(\phi_{w'})(I_{2n+1})\\
&=\il{N_w\bs N}\phi_{w'}(w^{-1}n)\,dn.
\endalign$$
This integral is zero unless $w^{-1}n\in B{w'}^{-1}\B$
for some $n\in N_w\bs N$.  If $w'=w_0$, it
is easy to see that this may happen only for $w=w_0$
and $n\in N\cap K$; part (1) follows.
Similarly, if $w'=w_1$, then $c_\chi(w,w')=0$ unless
$w=w_0$ or $w=w_1$, and $c_\chi(w_1,w_1)=1$.
Finally, to determine $c_\chi(w_1,w_0)$, it
suffices to take $g$ in (4.7) to be
of the form $\iota_n(g')$, with $g'\in GL(2,F)$.
In that case the $\phi_{w_i}$ and $f_{w_i}$ necessarily
match the analogously defined $GL(2)$ functions
in $\Ind(\chi_n,\chi_n^{-1})$,
and the determination of this coefficient follows
from the analogous statement on $GL(2,F)$, which
is equivalent to (2.5).

To prove the remaining parts of the Lemma,
it suffices to show that
$$L(\phi_{w_0})=0,\qquad L(\phi_{w_1})=1.\tag4.8$$
However, for $w\in\Omega$,
$$L(\phi_w)=\il{U}\phi_w(w_1u)\,\theta_S(u)^{-1}\,du.$$
A matrix calculation shows that
$w_1u\in B w^{-1}\B$ is impossible if $w=w_0$,
and holds when $w=w_1$ if and only if
$u\in U\cap K$.  Thus (4.8) holds.
\smudge
\enddemo

Combining (4.6), Lemma 4.1, and the evaluation of $R(d,w;\chi)$,
one finds that
$$H_\chi(d)=
c_{w_1}(\chi)\,\sigma_{^{w_1}\!\chi}(d)
+ \sum_{w\in \Omega;w\neq w_0,w_1}R(d,w;\chi)\,L(f_w).
\tag4.9$$
The final value for $\H_\chi(d)$ may now be
obtained by taking the first term on the right hand
side of (4.9), multiplying by the normalizing factor in the
functional equation of Theorem 1.3, and symmetrizing with respect
to the action of $\Omega$ on the parameters $\alpha_i$.
To do this, observe that (4.4) gives
$$c_{w_1}(\chi)=\prod_{1\leq i<j\leq n}
\frac{(1-\alpha_i\alpha_jq^{-1})
(1-\alpha_i\alpha_j^{-1}q^{-1})}
{(1-\alpha_i\alpha_j)
(1-\alpha_i\alpha_j^{-1})}
\prod_{1\leq i< n}
\frac{1-\alpha_i^2q^{-1}}{1-\alpha_i^2}.\tag4.10$$
Substituting this formula and (4.5) into (4.9)
and making use of Weyl's identity (1.17), the explicit formula follows.

This completes the proof of Theorem 1.5. \smudge
\smallskip
We turn to the proof of Theorem 1.6.  Suppose that $T$ is split,
so that the Bessel functional $B$ is given by equation (1.14).
Recall that $N$ denotes the unipotent radical of the standard
Borel subgroup of $G$.  Let $N(\O)=N\cap K$.
Given $k=(k_1,\cdots,k_n)$ with all $k_i\geq0$, define
$$P_k(g)=\il{N(\O)}\Phi_\chi(g n d_k)\,dn.$$
Then using the bijection (4.3),
an argument similar to the one given there demonstrates that
$P_k\in\Ind(\chi)^\B$.

Though we will ultimately use the Casselman-Shalika method,
we first establish the following Lemma, which implies
that the determination of $h(k_1,\cdots,k_n)$ follows from the
determination of the quantities $B(P_k)$.  For convenience,
let us set
$$H(k):=H_\chi(g_k)\qquad B(k):=B_\chi(P_k).$$
Then we have
\proclaim{Lemma~4.2} Let $k=(k_1,\cdots,k_n)$ with all $k_i\geq0$.
\roster
\item Suppose $k_1=0$.  Then $H(k)=B(k)$.
\item Suppose $k_1>0$.  Then
$$H(k)=(1-q^{-1})^{-1}\,(B(k)-q^{-1}\beta
B(k_1-1,k_2+1,k_3,\cdots,k_n)).$$
\endroster
\endproclaim

\demo{Proof} Comparing the definitions,
$$B(k)=\il{\O}B\left(\pi\left(n(x)d_k\right)\Phi_\chi\right)\,dx.$$
A matrix calculation shows that if $x\in\O$, then
$$n(x)d_k=\cases d_k\,n(\varpi^{-k_1}x)&\text{if $x\in\varpi^{k_1}\O$}\\
t(x)\,n(1)\,d_{k_1-m,k_2+m,k_3,\cdots,k_m}\,t(\varpi^m x^{-1})&
\text{if $x\in\varpi^m\O^\times, 0\leq m<k_1$.}
\endcases$$
If $x\in\varpi^{k_1}\O$ then $n(\varpi^{-k_1}x)\in K$;
moreover $\Phi_\chi$ is $K$-fixed.  Factoring
$$d_k=t(\varpi^{k_1})\,d_{0,k_1+k_2,k_3,\cdots,k_m}$$
and using property (1.10),
one thus obtains
$$\il{\varpi^{k_1}\O}B\left(\pi\left(n(x)d_k\right)\Phi_\chi\right)\,dx
=q^{-k_1}\beta^{k_1}H(0,k_1+k_2,k_3,\cdots,k_n).$$
Similarly if $x\in\varpi^m\O^\times$, then $t(\varpi^m x^{-1})\in K$.
Arguing similarly, one then finds that
$$\multline
B(k)=q^{-k_1}\beta^{k_1}H(0,k_1+k_2,k_3,\cdots,k_n)\\
+(1-q^{-1})\sum_{m=0}^{k_1-1}q^{-m}\beta^m H(k_1-m,k_2+m,k_3,\cdots,k_n).
\endmultline
\tag4.11$$
If $k_1=0$, Lemma 4.2, part (1), follows at once from (4.11).
For $k_1>0$, equation (4.11) also implies that
$$\multline q^{-1}\beta B(k_1-1,k_2+1,k_3,\cdots,k_n)=
q^{-k_1}\beta^{k_1}H(0,k_1+k_2,k_3,\cdots,k_n)\\
+(1-q^{-1})\sum_{m=1}^{k_1-1}q^{-m}\beta^m H(k_1-m,k_2+m,k_3,\cdots,k_n).
\endmultline
\tag4.12$$
Subtracting (4.12) from (4.11), part (2) of the Lemma follows.
\smudge
\enddemo

Write now $P_k$ in terms of the Casselman basis $\{f_w\}$:
$$P_k=\sum_{w\in\Omega}S(k,w;\chi)\,f_w.$$
The coefficients $S(k,w;\chi)$ are given by
$$\align
S(k,w;\chi)&=\il{N_w\bs N}P_k(w^{-1}n)\,dn\\
&=\il{N_w\bs N}\Phi_\chi(w^{-1}nd_k)\,dn\\
&=c_w(\chi)\,(^w\!\chi \delta_B^{1/2})(d_k),\tag4.13
\endalign$$
where $c_w(\chi)$ is defined by equation (4.4) above.

In this case, unlike the nonsplit case, no two terms for different $w$
contribute the same rational function of the $\alpha_i$.  Hence it
suffices to determine $B(f_{w_0})$; the value of $B(P_k)$ is then
obtained by symmetrization, using the functional equation of Theorem
1.4.  The value of $B(f_{w_0})$ is given by

\proclaim{Lemma~4.3} Suppose $|\alpha_n\beta|<q^{1/2}$.  Then
$$B(f_{w_0})=\left(1-\alpha_n\beta q^{-1/2}\right)^{-1}.$$
\endproclaim

\demo{Proof}  Suppose $|\alpha_n\beta|<q^{1/2}$.
Since $T$ normalizes $U$ and fixes $\theta_S$,
$$B(f_{w_0})=\il{F^\times}\il{U}f_{w_0}(w_0\,n(1)\,t(a)\,u)\,
\theta_S(u)^{-1}\,\lambda^{-1}(a)\,du\,d^\times a.$$
Suppose that $w_0n(1)t(a)u\in Bw_0\B$.
Then $(t(a^{-1})n(1)t(a))\,u\in w_0 B w_0\B$.
This relation implies that $u\in U\cap K$.  Thus
$$B(f_{w_0})=\il{F^\times}f_{w_0}(w_0\,n(1)\,t(a))\,
\,\lambda^{-1}(a)\,d^\times a.$$
Suppose first that $a\in F^\times$ satisfies $|a|_F<1$.  Then
one sees by a matrix calculation that in fact
$w_0n(1)t(a)$ is never is the cell $Bw_0\B$.  If
instead $|a|_F\geq1$, then
$$w_0\,n(1)\,t(a)=w_0\,t(a)\,(t(a)^{-1}n(1)t(a)),$$
and the last factor $t(a)^{-1}n(1)t(a)$ is in $\B$.
Since $f_{w_0}$ is right $\B$-invariant, this gives
$$\align
f_{w_0}(w_0n(1)t(a))&=f_{w_0}(w_0t(a)w_0^{-1}\cdot w_0)\\
&=\chi \delta_B^{1/2}\left(t(a^{-1})\right).
\endalign$$
Hence
$$\align
B(f_{w_0})&=\sum_{m=0}^\infty \alpha_n^m\,\beta^m\,q^{-m/2}\\
&=\left(1-\alpha_n\beta q^{-1/2}\right)^{-1},
\endalign
$$
as claimed.
\smudge
\enddemo

To conclude the proof of Theorem 1.6, observe that
$$B(k)=B_\chi(P_k)=\il{\O}H_\chi(n(x)d_k)\,dx.$$
By Theorem 1.4, this satisfies a functional equation under
the action of $\Omega$ on the $\alpha_i$.  Similarly
to (4.10), (4.4) gives
$$c_{w_0}(\chi)=\prod_{1\leq i<j\leq n}
\frac{(1-\alpha_i\alpha_jq^{-1})
(1-\alpha_i\alpha_j^{-1}q^{-1})}
{(1-\alpha_i\alpha_j)
(1-\alpha_i\alpha_j^{-1})}
\prod_{1\leq i\leq n}
\frac{1-\alpha_i^2q^{-1}}{1-\alpha_i^2}.$$
Substituting this expression into (4.13) to
obtain $S(k,w_0;\chi)$, and
symmetrizing, one obtains the pleasant formula
$$\multline
{{\prod_{i=1}^n(1-\alpha_i\beta q^{-1/2})(1-\alpha_i\beta^{-1} q^{-1/2})}
\over{
\prod_{1\leq i<j\leq n}(1-\alpha_i\alpha_jq^{-1})
(1-\alpha_i\alpha_j^{-1}q^{-1})\,\prod_{i=1}^n(1-\alpha_i^2 q^{-1})}}B(k)\\
=q^{e_k}\Delta^{-1}\times\\
\A\!\left((1-\alpha_n\beta^{-1}q^{-1/2})\alpha_1^{-k_n'-n}\prod_{i=1}^{n-1}
\alpha_{n+1-i}^{-k'_i-i}(1-\alpha_i\beta q^{-1/2})
(1-\alpha_i\beta^{-1}q^{-1/2}) \right).
\endmultline$$
Then applying Lemma 4.2, one obtains without difficulty the
expression for $h(k)$ given in Theorem 1.6.
(Note that if $k_1=0$ then in fact one obtains the formula
$$\multline h(0,k_2,\cdots,k_{n})=
q^{e_k}\,\Delta^{-1}\\\times
\A\!\left(\alpha_n^{-1}(1-\alpha_n\beta^{-1}q^{-1/2})
\prod_{i=1}^{n-1}
\alpha_{i}^{-k'_{n+1-i}-n-1+i}(1-\alpha_i\beta q^{-1/2})
(1-\alpha_i\beta^{-1}q^{-1/2})
\right).\endmultline$$
However, expanding the second factor, the term
$$\A\!\left(\alpha_n^{-1}(-\alpha_n\beta^{-1}q^{-1/2})
\prod_{i=1}^{n-1}
\alpha_{i}^{-k'_{n+1-i}-n-1+i}(1-\alpha_i\beta q^{-1/2})
(1-\alpha_i\beta^{-1}q^{-1/2})
\right)$$
is independent of $\alpha_n$, and hence its alternator
is zero. Thus one is in fact left with the expression
of Theorem 1.6, part (1).)

This completes the proof of Theorem 1.6.\smudge

\subheading{5.  Continuation of the Bessel Functionals:  The
Application of Bernstein's Theorem} Bernstein~\cite{Be} gave a powerful
new method for the meromorphic continuation of functionals satisfying a
suitable uniqueness property; in the case at hand, this uniqueness is
the uniqueness of the Bessel model, which, as we have already noted, was
proved by Novodvorsky~\cite{No}. Bernstein's result will appear as an
appendix to a book in progress of Cogdell and Piatetski-Shapiro; in
the meantime, a statement (but no proof) may be found in Gelbart and
Piatetski-Shapiro~\cite{GP}. In this Section we
will show that Bernstein's theorem
implies Theorem~1.7. We assume familiarity with either~\cite{Be} or
its paraphrase in~\cite{GP}.

Let $X=C^\infty\big((B\cap K)\bs K\big)$.
Given the $\alpha_i$, we may identify $X$ with
$V_\pi$ by extending an element $\Psi\in X$ uniquely to an element
of $\Psi_\chi\in V_\pi$ satisfying (1.11).
Let $D=(\C^\times)^{n+1}$. If
$(\alpha_1,\cdots,\alpha_n,\beta)\in D$, we consider functionals $B$
which satisfy the following two conditions:
\roster
\item $B$ is a Bessel functional,
\item $B(\Phi_\chi)=1.$
\endroster
We claim that these conditions may be expressed by a
polynomial system of equations (in Bernstein's sense) in the
parameters
$(\alpha_1,\cdots,\alpha_n,\beta)$. To see this, note that if $\Psi\in
X$ and (with the notation as in (1.10)) if $t\in T$, $u\in U$, then there
exist $\Psi_j\in X$ and polynomial functions $f_j$ ($j=1,\cdots,N$) of
the $\alpha_i$ and $\beta$ such that $\pi(tu)\Psi_\chi=\sum f_j\Psi_{j,\chi}$.
Thus (1.10) may be expressed by the polynomial equation
$$\sum_{j=1}^N f_jB(\Psi_j)-\theta_S(tu)B(\Psi)=0.$$
It is also evident that condition (2) above is a
polynomial equation. By Novodvorsky's theorem, the solution, if it
exists, is unique, and we have proved that a solution exists on a
non-empty open subset of $D$.
Consequently Bernstein's theorem is applicable, yielding the
meromorphic continuation of the functional $B$ in the sense made
precise by Theorem~1.7.\smudge

\subheading{6. Bessel Periods of Eisenstein Series}
In this Section we present a global application of these formulas.
(In fact, the application makes use of only a particular case of
them, namely the formula (6.6) for the value of the local
Bessel functional at the identity.)
Accordingly, we now let $F$ be a global field and ${\Bbb A}$ be
its ring of adeles.
Let $\pi=\otimes \pi_v$ be a cuspidal automorphic representation of
$GL(n,{\Bbb A})$.  Let $P$ be the standard maximal parabolic subgroup of
$G=SO(2n+1)$ with Levi factor $GL(n)$. Denoting by $\delta_P$ the
modular character of $P$, let
$f_s\in\Ind_{P_{\Bbb A}}^{G({\Bbb A})}(\pi\otimes\delta_P^{s-1/2})$, and let
$$E(g,s,f_s)=L_S(\pi,2n(s-1/2)+1,\vee^2)\,
\sum_{\gamma\in P_F\bs G_F}f_s(\gamma g)\tag6.1$$
be the Eisenstein series attached to $f_s$.
Here $S$ is a finite set of places including the archimedean ones and
those where $\pi$ is ramified, and $L_S(\pi,s,\vee^2)$ is the partial
symmetric square L-function, which is the normalizing factor of the
Eisenstein series. Let $a$, $b$ and $c$ be elements of $F$ such that
$b^2+2ac$ is not a square. Let $Q$ be the standard parabolic subgroup
of $G$ with Levi factor $GL(1)\times\cdots\times GL(1)\times SO(3)$,
let $U$ be the unipotent radical of $Q$, and let $\psi$ be a
nontrivial character of ${\Bbb A}/F$. Let $\theta$ be the character of
$U({\Bbb A})$
defined by (1.9). Then $Q({\Bbb A})$ acts on $U({\Bbb A})$ and hence on its
character group by conjugation; let $R({\Bbb A})$ be the subgroup
of the stabilizer consisting of elements whose projection to
the Levi factor of
$Q({\Bbb A})$ lies in the embedded $SO(3,{\Bbb A})$. We may
naturally extend $\theta$ to a character of $R({\Bbb A})$. Then
$R({\Bbb A})$ is the group of adelic points for an algebraic group $R$ which
is the semidirect product of $U$ with a one-dimensional torus $T$,
which for simplicity we are assuming nonsplit---this is our hypothesis
that $b^2+2ac$ is a nonsquare. In this case, there exists a unique
quadratic field extension $K$ of $F$ over which $T$ splits. Let
$\eta=\otimes \eta_v$
be the quadratic Hecke character of $F$ attached to $K$. In
this Section we will prove

\proclaim{Theorem 6.1}
The integral
$$\int_{R(F)\bs R({\Bbb A})}E(r,s,f_s)\,\theta(r)\,dr\tag6.2$$
is an Euler product, whose local factor at a good place $v$ equals
$$L\big(n(s-1/2)+1/2,\pi_v\big)\,L\big(n(s-1/2)+1/2,\pi_v\otimes\eta_v\big).
\tag6.3$$
\endproclaim

If $n=2$ this result is essentially due to B\"ocherer~\cite{B\"o} and
Mizumoto~\cite{Mi} (these authors consider holomorphic Siegel modular
forms on $PGSp_4$ over $\Q$, but the unramified local computation is
the same for general base field and infinity type).
The precise conditions to make $v$ good
are described below.

We turn to the proof of Theorem 6.1.
Unfolding the integral, we see that (6.2) equals
$$L_S(\pi,2n(s-1/2)+1,\vee^2)\,\sum_{\gamma\in P_F\bs G_F/R_F}
\int_{R^\gamma_{\Bbb A}\bs R_{\Bbb A}}\int_{R^\gamma_F\bs R^\gamma_{\Bbb A}}
f_s(\gamma ur)\,\theta(ur)\,du\,dr,$$
where $R^\gamma$ denotes the
algebraic group $R\cap\gamma^{-1}P\gamma$.  Using the fact that $\pi$
is cuspidal, one may show that only the open orbit in $P\bs G/R$
contributes, and we may take the representative
$$\gamma=\pmatrix&&I_n\\&(-1)^n\\I_n\endpmatrix.$$
Then $\gamma R^\gamma\gamma^{-1}=U$ is a maximal unipotent subgroup in
$GL(n)$, the Levi factor of $P$, and the character $\gamma
u\gamma^{-1}\mapsto\theta(u)$ is nondegenerate; indeed, since $b^2+2ac$
is a nonsquare, $a\ne0$. It follows that
$$W_s(g)=\int_{U_F\bs U_{\Bbb A}}f_s(ug)\,\theta(\gamma^{-1}u\gamma)\,du$$
lies in $\Ind_{P_{\Bbb A}}^{G_{\Bbb A}}
(\Cal W_\pi\otimes\delta_P^{s-1/2})$, where
$\Cal W_\pi$ is the Whittaker model of $\pi$ (relative to the appropriate
character of its maximal unipotent group). Writing
$$\Ind_{P_{\Bbb A}}^{G_{\Bbb A}}(\Cal W_\pi\otimes\delta_P^{s-1/2})=
\bigotimes_v \Ind_{P(F_v)}^{G(F_v)}(\Cal W_{\pi,v}\otimes\delta_P^{s-1/2})$$
as a restricted tensor product over all places $v$ of $F$, where
$\Cal W_{\pi,v}$ is the local Whittaker model of $\pi_v$, there is no
loss of generality in assuming that $W_s(g)$ is a pure tensor; thus we
write $W_s(g)=\prod_v W_{s,v}(g_v)$. The integral (6.2) thus equals
$$L_S(\pi,2n(s-1/2)+1,\vee^2)\int_{R^\gamma_{\Bbb A}\bs R_{\Bbb A}}
W_s(\gamma r)\,\theta(r)\,dr$$
and hence
is factorizable, with local factor
$$L(\pi_v,2n(s-1/2)+1,\vee^2)\int_{R^\gamma(F_v)\bs R(F_v)}
W_{s,v}(\gamma r)\,\theta(r)\,dr.\tag6.4$$

We compute this local factor at a good place.
More precisely, suppose that the finite place $v$ does
not ramify in $K$.  Then we compute (6.4) for any nonramified
principal series $\pi_v$
with Satake parameters $\alpha_1,\cdots,\alpha_n$,
under the assumption that $W_{s,v}$ is the
unramified spherical vector in
$\Ind_{P(F_v)}^{G(F_v)}(\Cal W_{\pi,v}\otimes\delta_P^{s-1/2})$,
normalized so that $W_{s,v}(1)=1$.
To carry out this computation we first make the further assumption that
$|\alpha_i|<|\alpha_{i+1}|$; the general case then follows
by analytic continuation.
(Note that the assumption $|\alpha_i|<|\alpha_{i+1}|$
is unrealistic for a representation $\pi_v$ which is a local component of an
automorphic representation $\pi$, as this condition violates
the Ramanujan conjecture.)
According to the results of Casselman and
Shalika~\cite{CS}, we may then write
$$W_{s,v}(g)=\prod_{1\leq i<j\leq n}\,(1-\alpha_i\alpha_j^{-1}q_v^{-1})^{-1}
\int_{U(F_v)}\Phi_{\chi,v}\left(\pmatrix J_n\\&1\\&&J_n\endpmatrix ug\right)
\theta(\gamma^{-1}u\gamma)\,du,$$
where $q_v$ is the cardinality of the residue field at $v$, $J_n$ is the
$n\times n$ matrix with ones on the sinister diagonal and zeros elsewhere,
and
$\Phi_{\chi,v}$ is the normalized spherical vector in the
representation denoted (in the notation of Section~1)
$\Ind(\chi_1,\cdots,\chi_n)$, where
$\chi_i(\varpi_i)=\alpha_i\,q^{-n(s-1/2)}$.  Substituting this into
(6.4), we obtain precisely the integral (1.12) if the place $v$ is
inert in $K$, that is, if $T$ is a nonsplit torus in $F_v$.
If on the other hand the place $v$ splits in $K$, we obtain the integral
(1.14). To see this, one must remember to conjugate the torus so as to
make it diagonal; and conjugating $w_1$ in this way produces $w_0n(1)$
as in (1.14). One sees that, in either case, the local
factor (6.4) is equal to
$$\prod_{1\leq i\le j\leq n}\,
\left(1-\alpha_i\alpha_j\,q_v^{-1-2n(s-1/2)}\right)^{-1}
\prod_{1\leq i<j\leq n}\,
\left(1-\alpha_i\alpha_j^{-1}q_v^{-1}\right)^{-1}\,H_\chi(1).
\tag6.5$$
By Theorems~1.5 and~1.6, $\H_\chi(1)=1$, so by Theorems~1.3 and~1.4
(with $\alpha_i\,q_v^{-n(s-1/2)}$ replacing $\alpha_i$, and
$\beta=1$),
$$\multline
H_\chi(1)=\\
\prod_{1\leq i<j\leq n}\,
\left(1-\alpha_i\alpha_j\,q_v^{-1-2n(s-1/2)}\right)\,
\left(1-\alpha_i\alpha_j^{-1}q_v^{-1}\right)
\prod_{i=1}^n\,
{1+\eta_v(\varpi_v)\alpha_i\,q_v^{-n(s-1/2)-1/2}\over
1-\alpha_i\,q_v^{-n(s-1/2)-1/2}},
\endmultline\tag6.6$$
where the quadratic character $\eta_v(\varpi_v)$ equals $1$ if $v$
splits in $K$, and
$-1$ if $v$ is inert.  Substituting this formula into (6.5) and simplifying,
we obtain the local factor (6.3) for $v$.\smudge

\newpage
\enddocument